\documentclass[sts,preprint]{imsart}
\usepackage{graphicx}	
\usepackage{natbib} 
\usepackage{amssymb, amsmath}
\usepackage[colorlinks,citecolor=blue,urlcolor=blue]{hyperref}

\usepackage{tikz}
\usetikzlibrary{decorations.pathreplacing, arrows.meta, math}

\definecolor{light}{RGB}{220, 188, 188}
\definecolor{mid}{RGB}{185, 124, 124}
\definecolor{dark}{RGB}{143, 39, 39}
\definecolor{highlight}{RGB}{180, 31, 180}
\definecolor{gray10}{gray}{0.1}
\definecolor{gray20}{gray}{0.2}
\definecolor{gray30}{gray}{0.3}
\definecolor{gray40}{gray}{0.4}
\definecolor{gray60}{gray}{0.6}
\definecolor{gray70}{gray}{0.7}
\definecolor{gray80}{gray}{0.8}
\definecolor{gray90}{gray}{0.9}
\definecolor{gray60}{gray}{0.95}

\tikzmath{
  function normal(\x) {
    return exp(-0.5 * \x * \x);
  };
  function laplace(\x) {
    return exp(-sqrt(\x * \x));
  };
}

\begin{document}

\begin{frontmatter}

\title{A Short Review of Ergodicity and Convergence of Markov chain Monte Carlo Estimators}
\runtitle{Ergodicity and Convergence}

\begin{aug}
  \author{Michael Betancourt \ead[label=e1]{betan@symplectomorphic.com}}

  \runauthor{M. Betancourt}

  \address{Michael Betancourt is the principal research scientist of Symplectomorphic, LLC. \printead{e1}.}

\end{aug}

\begin{abstract}
This short note reviews the basic theory for quantifying both the asymptotic and preasymptotic 
convergence of Markov chain Monte Carlo estimators.
\end{abstract}

\end{frontmatter}

\maketitle

In this note we'll review useful notions of distance between probability distributions 
and use them to specify the conditions under which a Markov transition will converge 
to a unique invariant distribution, both asymptotically and preasymptotically.  The 
material will assume familiarity with the basic concepts of measure theory.

\section{Markov, Markov, Markov!}

Consider an ambient topological space $X$ equipped with the $\sigma$-algebra 
$\mathcal{X}$.  A \emph{Markov transition}
\begin{alignat*}{6}
\tau :\; &\mathcal{X} \times X& &\rightarrow& \; &\mathbb{R}&
\\
&A, x& &\mapsto& &\tau(A \mid x)&
\end{alignat*}
specifies a probability distribution at each point $x \in X$, and sampling 
from this distribution realizes a transition from that initial point $x_{0}$ to 
a new point $x_{1}$.  Applying the Markov transition repeatedly yields a sequence
of points $(x_{0}, x_{1}, \ldots, x_{N})$ called a \emph{Markov chain} that meanders 
through the ambient space.

More formally a Markov transition lifts any distribution of initial points $\rho$
into a joint distribution over the product space of initial points and transitions 
$(x_{0}, x_{1}) \in X \times X$, $\tau \times \rho$.  Applying the Markov transition 
again yields a joint distribution over the product space of three states 
$(x_{0}, x_{1}, x_{2}) \in X \times X \times X$,
\begin{equation*}
\tau^{2} \times \rho = \tau \times (\tau \times \rho) = \tau \times \tau \times \rho,
\end{equation*}
and applying it $N$ times yields a joint distribution over possible Markov chain 
configurations $(x_{0}, \ldots x_{N}) \in X^{N + 1}$,
\begin{equation*}
\tau^{N} \times \rho = \underbrace{\tau \times \ldots \times \tau}_{N\text{ times}} \rho.
\end{equation*}

Lifting $\rho$ to the joint distribution $\tau \times \rho$ and then marginalizing 
to the space of transitions convolves the Markov transition over all possible 
initial states, yielding the \emph{$1$-step distribution} of possible transitions,
\begin{equation*}
\tau \circ \rho \equiv (\varpi_{1})_{*} (\tau \times \rho),
\end{equation*}
where $\varpi_{n}$ denotes the product space projection operator,
\begin{alignat*}{6}
\varpi_{n} :\; &X^{N + 1}& &\rightarrow& \; &X&
\\
&(x_{0}, \ldots, x_{n}, \ldots, x_{N})& &\mapsto& &x_{n}&.
\end{alignat*}

Repeating this convolution to the $1$-step distribution yields the marginal 
distribution of second transitions, or the $2$-step distribution,
\begin{equation*}
\tau^{2} \circ \rho = \tau \circ (\tau \circ \rho),
\end{equation*}
and repeating the convolution $N$ times yields the marginal distribution for
the $N + 1$ state in the Markov chain or the $N$-step distribution,
\begin{equation*}
\tau^{N} \circ \rho = \underbrace{\tau \circ \ldots \circ \tau}_{N\text{ times}} \, \rho.
\end{equation*}
Equivalently we can define the $N$-step distribution as the marginal distribution 
of the full joint distribution,
\begin{equation*}
\tau^{N} \circ \rho = (\varpi_{N})_{*} (\tau^{N} \times \rho).
\end{equation*}

A probability distribution $\pi$ that is invariant to the Markov convolution,
\begin{equation*}
\tau \circ \pi = \pi,
\end{equation*}
is said to be preserved by the Markov transition in which case $\pi$ is denoted 
a \emph{stationary} or \emph{invariant} distribution.  When a Markov transition 
is engineered to preserve a given distribution that distribution is also denoted 
the \emph{target} distribution.

We often take for granted that when a Markov transition $\tau$ preserves a given 
target distribution $\pi$ the Markov chains generated from $\tau$ will explore 
$\pi$.  Unfortunately invariance alone is not guaranteed to fully quantify $\pi$, 
even asymptotically as the chains grow to infinite lengths.  In order for Markov 
chains to adequately characterize $\pi$ the $N$-step distribution needs to converge 
towards $\pi$ sufficiently quickly.  To formally define the desired convergence
let alone speed of that convergence, however, we first need to define some notion 
of distance between probability distributions.

\section{Quantifying Convergence}

There are many notions of distance between measures, and hence probability 
distributions, but a few are particularly useful for the theoretical analysis of 
Markov chain Monte Carlo.  In this section we'll discuss integral probability 
metrics, Wasserstein metrics, their duality, and important special cases.

\subsection{Integral Probability Metrics}

An integral probability metric \cite{Muller:1997} defines a notion of distance 
between measures by differences in the expectation values of certain functions.

Consider two measures $\mu$ and $\nu$ over $X$ and a space of real-valued, 
measurable functions,  $f : X \rightarrow \mathbb{R} \in \mathcal{F}$.  If both
$\mathbb{E}_{\mu}[f]$ and $\mathbb{E}_{\nu}[f]$ are finite for all 
$f \in \mathcal{F}$ then we can define the integral probability metric for 
$\mathcal{F}$ as the largest difference between expectation values,
\begin{equation*}
|| \mu - \nu ||_{\mathcal{F}}
= \underset{f \in \mathcal{F}}{\mathrm{sup}} 
\left| \mathbb{E}_{\mu}[f] - \mathbb{E}_{\nu}[f] \right|.
\end{equation*}
One immediate consequence of this definition is that an integral probability 
metric bounds how well expectations with respect to $\mu$ approximate expectations 
with respect to $\nu$, and vice versa, within the space of test functions $\mathcal{F}$.

When $\mathcal{F}$ spans functions of practical interest any constraints on the 
corresponding integral probability metric are directly applicable in applied settings.  
More restrictive spaces, however, can still prove useful more indirectly.  In particular 
integral probability metric bounds for restricted spaces of functions can sometimes
demonstrate the absence of pathological behavior that would otherwise effect the 
expectation values of functions outside of that space.

\subsection{Wasserstein Metric} \label{sec:wasserstein}

Wasserstein metrics \cite{SriperumbudurEtAl:2009} define a notion of distance between 
probability distributions through a distance function defined on the ambient space.  
In particular let $g : X \times X \rightarrow \mathbb{R}^{+}$ be a positive, symmetric 
function, $g(x_{1}, x_{2}) = g(x_{2}, x_{1}) \ge 0$, that vanishes if and only if the 
two arguments are identical, $g(x_{1}, x_{2}) = 0$ if and only if $x_{1} = x_{2}$.

We first assume that the distance function is compatible with the two measures being
compared so that the expectation values of both arguments are well-defined,
\begin{align*}
\mathbb{E}_{\mu}[g(-, x)] &= \mathbb{E}_{\mu}[g(x, -)] < \infty, \, \forall x \in X
\\
\mathbb{E}_{\nu}[g(-, x)] &= \mathbb{E}_{\nu}[g(x, -)] < \infty, \, \forall x \in X.
\end{align*}

In order to take an expectation value of both arguments of the distance function we 
need to lift $\mu$ and $\nu$ together to a joint measure on the product space $X \times X$.  
A \emph{coupling} of $\mu$ and $\nu$ is any joint distribution $\gamma$ over $X \times X$ 
that marginalizes to $\mu$ and $\nu$ in each component,
\begin{align*}
\gamma(A \times X) &= \mu(A), \forall A \in \mathcal{X}
\\
\gamma(X \times A) &= \nu(A), \forall A \in \mathcal{X}.
\end{align*}
We'll denote the space of all couplings of $\mu$ and $\nu$ by  $\Gamma(\mu, \nu)$.
Because we assumed that the marginal expectation values are well-defined the joint
expectation value with respect to any coupling $\mathbb{E}_{\gamma}[g] \in \mathbb{R}$
will also be well-defined.

The \emph{$1$-Wasserstein metric}, also known as the \emph{Kantorovich-Rubinstein distance}, 
is then defined as the smallest of these joint expectation values over all possible couplings,
\begin{equation*}
W_{1, g}(\mu, \nu) = \underset{ \gamma \in \Gamma(\mu, \nu) }{ \mathrm{inf} }
\mathbb{E}_{\gamma}[g].
\end{equation*}

\subsection{Kantorovich-Rubinstein Theorem}

Although they might appear different superficially, integral probability metrics 
and $1$-Wasserstein metrics actually provide equivalent notions of distance between 
two measures.  This duality is formalized in the \emph{Kantorovich-Rubinstein Theorem}
\cite{KantorovicEtAl:1958, Dudley:2002}.

A distance function $g$, defines a \emph{Lipschitz semi-norm} on functions by
\begin{equation*}
|| f ||_{g} 
= \underset{ x_{1} \ne x_{2} }{ \mathrm{sup} } 
\frac{ | f(x_{1}) - f(x_{2}) | }{ g(x_{1}, x_{2}) }.
\end{equation*}
The Lipschitz semi-norm in turn defines a space of all continuous functions whose 
semi-norms are less than or equal to $1$,
\begin{equation*}
\mathcal{F}_{g} = \left\{ f \in C^{0}(X) \mid \, || f ||_{g}  \le 1 \right\}.
\end{equation*}

While the positional metric immediately defines a 1-Wasserstein metric $W_{1, g}(\mu, \nu)$, 
this space of Lipschitz functions defines a corresponding integral probability metric, 
$|| \mu - \nu ||_{\mathcal{F}_{g}}$.  The Kantorovich-Rubinstein Theorem shows that 
these two metrics are in fact equal provided that all of the relevant expectation values 
are well-defined.

The ability to jump back and forth between integral probability and 1-Wasserstein 
metrics is particularly useful in practice.  For example one might be able to 
theoretically bound convergence in one metric but then use its dual representation 
to understand the consequences of that bound in practical settings.

\subsection{The Total Variation Distance}

The total variation distance is a particularly useful metric for the analysis of Markov 
chain Monte Carlo that can be derived from either of these dual perspectives. 

From the 1-Wasserstein perspective the total variation distance is given by 
using a degenerate distance function specified by an indicator function,
\begin{equation*}
g = \mathbb{I}_{D^{c}},
\end{equation*}
where
\begin{equation*}
D = \left\{ (x_{1}, x_{2}) \in X \times X \mid x_{1} = x_{2} \right\}. 
\end{equation*}
and
\begin{equation*}
D^{c} = \left\{ (x_{1}, x_{2}) \in X \times X \mid x_{1} \ne x_{2} \right\}.
\end{equation*}
In other words the distance between any two distinct points is always 1 while 
the distance between any point and itself is 0.  Working through the coupling 
definition of the $1$-Wasserstein metric one can show that the total variation
distance is equal to the largest difference of allocated measure over elements 
of the defining $\sigma$-algebra,
\begin{equation*}
\mathrm{TV}(\mu, \nu) 
= || \mu - \nu  ||_{\mathrm{TV}} 
= \underset{ A \in \mathcal{X} }{ \mathrm{sup} } | \mu(A) - \nu(A) |.
\end{equation*}

The corresponding integral probability metric definition is given by
\begin{equation*}
|| \mu - \nu  ||_{\mathrm{TV}} 
= \underset{ 0 \le f \le 1 }{ \mathrm{sup} } 
\left| \mathbb{E}_{\mu}[f] - \mathbb{E}_{\nu}[f] \right|,
\end{equation*}
where $0 \le f \le 1$ denotes the space of continuous, real-valued functions 
whose outputs are bounded between 0 and 1.

This bound on function behavior, and hence the class of test functions being
considered, can also be scaled so long as the normalization is modified.  More 
formally
\begin{equation*}
|| \mu - \nu  ||_{\mathrm{TV}} 
= \frac{1}{f_{\max} - f_{\min}}
\underset{ f_{\min} \le f \le f_{\max} }{ \mathrm{sup} } 
\left| \mathbb{E}_{\mu}[f] - \mathbb{E}_{\nu}[f] \right|.
\end{equation*}

Taking advantage of this scaling freedom the total variation distance is often 
defined as
\begin{align*}
|| \mu - \nu  ||_{\mathrm{TV}} 
&= \frac{1}{2}
\underset{ -1 \le f \le 1 }{ \mathrm{sup} } 
\left| \mathbb{E}_{\mu}[f] - \mathbb{E}_{\nu}[f] \right|
\\
&= \frac{1}{2}
\underset{ | f | \le 1 }{ \mathrm{sup} } 
\left| \mathbb{E}_{\mu}[f] - \mathbb{E}_{\nu}[f] \right|.
\end{align*}

\subsection{Application to Markov chain Monte Carlo}

Any of the probability metrics introduced above can provide some quantification
for how similar the $N$-step distribution of a Markov transition $\tau^{N} \circ \rho$
is to the target distribution $\pi$.  If the two distributions become more similar 
as $N$ increases then any Markov chain generated by the Markov transition will 
provide a better approximations to $\pi$ as it becomes longer.

The challenge at hand is then to show not only that the distance $|| \tau^{N} \circ \rho_{1} - \pi ||$
converges towards $0$ as $N$ increases but also provide an explicit bound on how 
quickly the distance decreases with increasing $N$.

\section{Asymptotic Convergence} \label{sec:asymptotic_convergence}

The invariance of the target distribution provides some constraints on how the 
$N$-step distribution of a Markov chain converges towards the target distribution
\cite{RobertsEtAl:2004}.  For example the total variation distance between 
$N$-step distributions from different initializations will be bounded by the distance 
between those initializations,
\begin{equation*}
|| \tau^{N} \circ \rho - \tau^{N} \circ \rho' ||_{\mathrm{TV}} \le || \rho - \rho' ||_{\mathrm{TV}}.
\end{equation*}
If we initialize the second Markov chain from stationarity, $\rho' = \pi$, then 
$\tau^{N} \circ \pi = \pi$ and this bound implies
\begin{equation*}
|| \tau^{N} \circ \rho - \pi ||_{\mathrm{TV}} \le || \rho - \pi ||_{\mathrm{TV}}.
\end{equation*}
In words the $N$-step distributions cannot move any further from the target 
distribution in total variation distance than the initial distribution.  

Moreover if we take $\rho = \tau^{M} \circ \rho'$ for some initialization $\rho'$ then 
\begin{equation*}
|| \tau^{N + M} \circ \rho' - \pi ||_{\mathrm{TV}} \le 
|| \tau^{M} \circ \rho' - \pi ||_{\mathrm{TV}}.
\end{equation*}
Consequently once the $N$-step distribution achieves a certain distance it can never
move further than that distance in future iterations.  In other words the $N$-step 
distributions are non-expanding in the total variation distance.

This non-expansive behavior of the total variational distance can also be interpreted 
as a \emph{data processing inequality} \cite{CoverEtAl:2006}.  Unfortunately 
non-expansive behavior in other probability metrics is not guaranteed by the invariance 
of the target distribution alone.

That said non-expansion in the total variation distance is still not enough to ensure 
that the total variation distance vanishes asymptotically.  In particular if the Markov 
transition features multiple invariant distributions then we cannot guarantee that 
Markov chains will focus on the desired target distribution, $\pi$.  The $N$-step 
distribution is guaranteed to converge to $\pi$ asymptotically only when certain 
undesired behavior can be avoided \cite{RobertsEtAl:2004}.

\subsection{Reducibility}

One important undesired behavior is that Markov chain initialized within some neighborhood 
might not be able to explore every other relevant neighborhood in $X$.  To formalize relevance 
here we define the $\pi$-null sets as those measurable neighborhoods with vanishing target 
probability,
\begin{equation*}
\mathcal{X}_{\pi\text{-null}} = \{ A \in \mathcal{X} \mid \pi(A) = 0 \},
\end{equation*}
and the complementary $\pi$-non-null sets as those measurable neighborhoods with non-zero
target probability
\begin{equation*}
\mathcal{X}_{\pi\text{-non-null}} = \{ A \in \mathcal{X} \mid \pi(A) > 0 \}.
\end{equation*}
In order to completely explore the target distribution every Markov chain needs to be able 
to eventually explore every $\pi$-non-null set.  

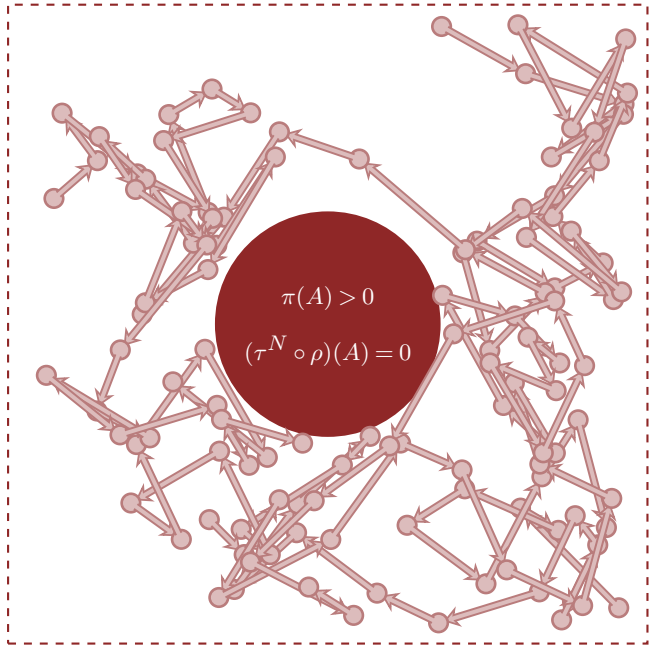
\begin{figure}
\centering
\begin{tikzpicture}[scale=1.0, thick]
  \draw[dark, dashed] (-4.25, -4.25) rectangle (4.25, 4.25);
 
  \fill[dark] (0, 0) circle (1.5);
  \node[white] at (0, 0.35) { $\pi(A) > 0$ };
  \node[white] at (0, -0.35) { $(\tau^{N} \circ \rho)(A) = 0$ };
 
  \pgfmathsetmacro{\N}{29}
 
  \foreach \x/\y [count=\n] in {1.513/3.960, 2.630/3.332, 3.938/2.926, 2.972/2.224, 3.931/2.792, 3.016/1.708, 1.984/1.104, 3.019/0.482, 3.661/0.814, 2.993/1.432, 3.609/2.174, 3.984/3.074, 2.528/3.982, 3.243/2.607, 3.961/3.794, 3.544/2.554, 2.641/1.154, 3.798/0.327, 3.236/1.237, 3.906/0.429, 2.586/1.539, 1.796/0.949, 2.154/-0.350, 2.890/0.418, 1.830/0.990, 0.419/2.197, -0.643/2.558, -1.381/1.423, -2.467/0.131} {
    \ifnum \n>1
      \fill[color=mid] (A) circle (4pt);
      \fill[color=light] (A) circle (3pt);
      \draw[color=mid, -{Stealth[length=6pt, width=8pt]}, 
            line width=3, shorten <=2.25pt, shorten >=2pt] (A) -- (\x, \y);
      \draw[color=light, -{Stealth[length=3pt, width=4pt]}, 
            line width=1.5, shorten >=3.25pt] (A) -- (\x, \y);
    \fi
    
    \ifnum \n=\N
      \fill[color=mid] (\x, \y) circle (4pt);
      \fill[color=light] (\x, \y) circle (3pt);
    \fi
    \coordinate (A) at (\x, \y);
  }
  
  \foreach \x/\y [count=\n] in {3.868/-3.774, 2.562/-2.448, 3.065/-2.744, 1.821/-2.191, 1.055/-2.676, 2.103/-3.455, 2.847/-2.031, 3.328/-1.273, 3.704/-2.706, 3.109/-3.935, 3.589/-2.937, 3.289/-2.549, 2.813/-3.562, 1.473/-3.974, 0.655/-3.579, -0.415/-2.775, -1.242/-3.530, -0.399/-2.394, 0.564/-1.493, 0.194/-1.875, -1.151/-2.730, -0.765/-3.067, -1.439/-1.690, -2.617/-2.380, -1.947/-2.863, -2.543/-1.602, -1.633/-0.335, -0.804/-1.772, -2.059/-0.768} {
    \ifnum \n>1
      \fill[color=mid] (A) circle (4pt);
      \fill[color=light] (A) circle (3pt);
      \draw[color=mid, -{Stealth[length=6pt, width=8pt]}, 
            line width=3, shorten <=2.25pt, shorten >=2pt] (A) -- (\x, \y);
      \draw[color=light, -{Stealth[length=3pt, width=4pt]}, 
            line width=1.5, shorten >=3.25pt] (A) -- (\x, \y);
    \fi
    
    \ifnum \n=\N
      \fill[color=mid] (\x, \y) circle (4pt);
      \fill[color=light] (\x, \y) circle (3pt);
    \fi
    \coordinate (A) at (\x, \y);
  }
  
    \foreach \x/\y [count=\n] in {-3.640/1.672, -3.062/2.175, -3.532/2.807, -2.531/1.999, -1.590/1.487, -2.113/2.788, -1.539/3.128, -1.025/2.811, -2.182/2.442, -1.531/1.408, -1.786/1.077, -2.419/1.928, -3.037/2.495, -2.555/1.792, -1.475/1.049, -0.697/2.224, -1.593/0.726, -2.425/0.284, -1.939/1.499, -1.608/1.053, -2.758/-0.343, -3.068/-1.156, -2.364/-1.513, -3.741/-0.681, -2.762/-1.473, -1.472/-1.081, -1.046/-1.878, -1.412/-1.268, -0.335/-1.583} {
    \ifnum \n>1
      \fill[color=mid] (A) circle (4pt);
      \fill[color=light] (A) circle (3pt);
      \draw[color=mid, -{Stealth[length=6pt, width=8pt]}, 
            line width=3, shorten <=2.25pt, shorten >=2pt] (A) -- (\x, \y);
      \draw[color=light, -{Stealth[length=3pt, width=4pt]}, 
            line width=1.5, shorten >=3.25pt] (A) -- (\x, \y);
    \fi
    
    \ifnum \n=\N
      \fill[color=mid] (\x, \y) circle (4pt);
      \fill[color=light] (\x, \y) circle (3pt);
    \fi
    \coordinate (A) at (\x, \y);
  }
  
    \foreach \x/\y [count=\n] in {-1.575/-2.601, -1.111/-3.066, -0.645/-2.337, -1.448/-3.637, 0.040/-2.865, 0.968/-1.582, 1.785/-1.940, 2.379/-3.269, 3.388/-3.744, 3.775/-2.319, 2.820/-1.862, 2.480/-0.783, 3.025/-1.714, 2.336/-1.082, 1.526/0.383, 2.489/0.109, 3.084/-0.517, 2.663/-0.173, 3.036/-0.885, 2.158/-0.517, 2.865/-1.712, 3.414/-0.552, 3.017/0.309, 1.673/-0.129, 0.823/-1.621, -0.181/-2.350, -1.027/-3.166, 0.344/-3.875, -0.254/-3.498} {
    \ifnum \n>1
      \fill[color=mid] (A) circle (4pt);
      \fill[color=light] (A) circle (3pt);
      \draw[color=mid, -{Stealth[length=6pt, width=8pt]}, 
            line width=3, shorten <=2.25pt, shorten >=2pt] (A) -- (\x, \y);
      \draw[color=light, -{Stealth[length=3pt, width=4pt]}, 
            line width=1.5, shorten >=3.25pt] (A) -- (\x, \y);
    \fi
    
    \ifnum \n=\N
      \fill[color=mid] (\x, \y) circle (4pt);
      \fill[color=light] (\x, \y) circle (3pt);
    \fi
    \coordinate (A) at (\x, \y);
  }
  
\end{tikzpicture}
\caption{
A $\pi$-reducible Markov transition $\tau$ initialized from $\rho$ features at least one set $A$ with 
non-zero invariant probability $\pi(A) > 0$ that cannot be reached by all Markov chains initialized 
from $\rho$ no matter how large $N$ is.
}
\label{fig:reducible_set} 
\end{figure}

\begin{figure}
\centering
\begin{tikzpicture}[scale=0.8, thick]

\begin{scope}[shift={(0, 0)}]
  \draw[dark, dashed] (-4.25, -4.25) rectangle (4.25, 4.25);
 
  \fill[dark] (-4.25, -4.25) rectangle (0, 4.25);
  \node[white] at (-2.125, 0.35) { $\pi(A_{-}) > 0$ };
  \node[white] at (-2.125, -0.35) { $(\tau^{N} \circ \rho_{+})(A_{-}) = 0$ };
 
  \pgfmathsetmacro{\N}{29}
 
  \foreach \x/\y [count=\n] in {1.868/2.205, 3.356/2.567, 2.938/1.831, 3.684/2.216, 3.040/0.975, 3.345/1.386, 1.931/2.370, 2.568/3.045, 1.185/1.801, 1.574/3.223, 2.704/3.649, 2.367/3.228, 3.420/2.582, 1.960/3.210, 0.505/2.348, 1.959/3.234, 1.642/2.467, 0.463/1.579, 1.417/0.357, 1.005/-0.284, 1.326/-0.623, 1.584/-1.817, 2.305/-3.120, 2.639/-1.831, 1.222/-3.188, 1.472/-2.670, 0.764/-1.668, 0.106/-0.242, 0.561/1.120} {
    \ifnum \n>1
      \fill[color=mid] (A) circle (4pt);
      \fill[color=light] (A) circle (3pt);
      \draw[color=mid, -{Stealth[length=6pt, width=8pt]}, 
            line width=3, shorten <=2.25pt, shorten >=2pt] (A) -- (\x, \y);
      \draw[color=light, -{Stealth[length=3pt, width=4pt]}, 
            line width=1.5, shorten >=3.25pt] (A) -- (\x, \y);
    \fi
    
    \ifnum \n=\N
      \fill[color=mid] (\x, \y) circle (4pt);
      \fill[color=light] (\x, \y) circle (3pt);
    \fi
    \coordinate (A) at (\x, \y);
  }
  
  \foreach \x/\y [count=\n] in {3.701/-3.525, 2.955/-2.655, 2.588/-3.525, 2.325/-2.722, 1.033/-1.712, 0.780/-3.197, 1.718/-3.623, 3.207/-3.969, 3.867/-3.674, 3.465/-2.381, 2.624/-0.888, 1.146/-0.569, 2.459/-1.179, 3.891/0.055, 2.698/0.375, 2.070/-0.693, 0.627/-2.026, 1.455/-3.320, 2.490/-2.371, 1.406/-3.435, 1.802/-3.959, 2.183/-3.361, 0.853/-2.393, 0.280/-3.118, 0.698/-2.154, 0.227/-1.288, 0.650/-2.079, 1.014/-3.399, 1.859/-2.251} {
    \ifnum \n>1
      \fill[color=mid] (A) circle (4pt);
      \fill[color=light] (A) circle (3pt);
      \draw[color=mid, -{Stealth[length=6pt, width=8pt]}, 
            line width=3, shorten <=2.25pt, shorten >=2pt] (A) -- (\x, \y);
      \draw[color=light, -{Stealth[length=3pt, width=4pt]}, 
            line width=1.5, shorten >=3.25pt] (A) -- (\x, \y);
    \fi
    
    \ifnum \n=\N
      \fill[color=mid] (\x, \y) circle (4pt);
      \fill[color=light] (\x, \y) circle (3pt);
    \fi
    \coordinate (A) at (\x, \y);
  }
\end{scope}

\begin{scope}[shift={(10, 0)}]
  \draw[dark, dashed] (-4.25, -4.25) rectangle (4.25, 4.25);
 
  \fill[dark] (0, -4.25) rectangle (4.25, 4.25);
  \node[white] at (2.125, 0.35) { $\pi(A_{+}) > 0$ };
  \node[white] at (2.125, -0.35) { $(\tau^{N} \circ \rho_{-})(A_{+}) = 0$ };
 
  \pgfmathsetmacro{\N}{29}
 
  \foreach \x/\y [count=\n] in {-2.564/2.665, -3.981/3.130, -2.652/2.798, -3.882/1.681, -2.678/1.394, -3.902/2.707, -3.605/3.338, -3.944/2.188, -2.613/3.648, -1.868/3.103, -2.695/1.878, -2.964/1.366, -3.337/-0.115, -1.909/0.704, -1.029/-0.248, -1.487/-0.614, -1.850/-1.107, -2.637/0.230, -3.383/-0.918, -2.561/0.556, -3.800/-0.671, -2.951/0.297, -1.456/1.764, -0.883/0.773, -1.707/-0.282, -2.495/-1.143, -1.882/-1.983, -1.139/-3.293, -0.158/-3.814} {
    \ifnum \n>1
      \fill[color=mid] (A) circle (4pt);
      \fill[color=light] (A) circle (3pt);
      \draw[color=mid, -{Stealth[length=6pt, width=8pt]}, 
            line width=3, shorten <=2.25pt, shorten >=2pt] (A) -- (\x, \y);
      \draw[color=light, -{Stealth[length=3pt, width=4pt]}, 
            line width=1.5, shorten >=3.25pt] (A) -- (\x, \y);
    \fi
    
    \ifnum \n=\N
      \fill[color=mid] (\x, \y) circle (4pt);
      \fill[color=light] (\x, \y) circle (3pt);
    \fi
    \coordinate (A) at (\x, \y);
  }
  
  \foreach \x/\y [count=\n] in {-3.959/-3.619, -3.036/-2.649, -2.419/-3.604, -1.655/-2.133, -1.157/-2.639, -2.062/-2.196, -2.571/-3.695, -1.118/-3.324, -2.334/-2.291, -1.286/-2.981, -0.984/-1.719, -0.483/-2.079, -1.523/-1.712, -2.621/-0.222, -1.823/1.159, -0.618/1.860, -1.208/0.538, -0.280/-0.621, -1.352/0.786, -0.030/1.443, -0.719/1.796, -1.383/1.482, -1.106/0.883, -2.026/0.025, -1.233/-1.102, -0.842/-2.233, -0.529/-3.731, -0.857/-3.271, -1.738/-3.761} {
    \ifnum \n>1
      \fill[color=mid] (A) circle (4pt);
      \fill[color=light] (A) circle (3pt);
      \draw[color=mid, -{Stealth[length=6pt, width=8pt]}, 
            line width=3, shorten <=2.25pt, shorten >=2pt] (A) -- (\x, \y);
      \draw[color=light, -{Stealth[length=3pt, width=4pt]}, 
            line width=1.5, shorten >=3.25pt] (A) -- (\x, \y);
    \fi
    
    \ifnum \n=\N
      \fill[color=mid] (\x, \y) circle (4pt);
      \fill[color=light] (\x, \y) circle (3pt);
    \fi
    \coordinate (A) at (\x, \y);
  }
\end{scope}
  
\end{tikzpicture}
\caption{
A $\pi$-reducible Markov transition $\tau$ often induces a partition of the ambient space.  
Here Markov chains initialized at positive values will never reach negative values while 
Markov chains initialized at negative values will never reach positive values.
}
\label{fig:reducible_partition} 
\end{figure}
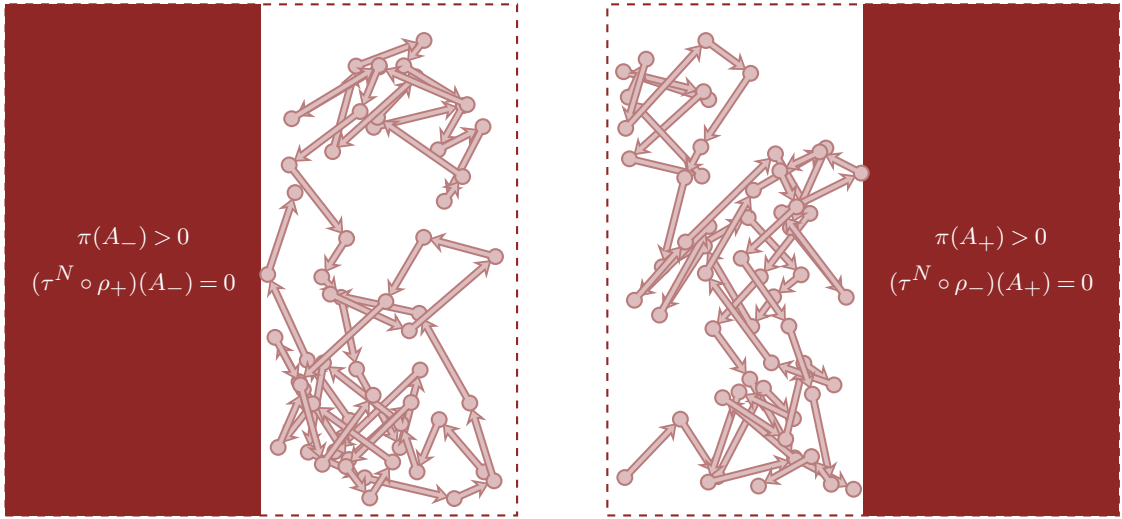

Markov transitions that generate Markov chains that cannot reach every relevant set, more 
formally when there is at least one $\pi$-non-null set $A \in \mathcal{X}_{\pi\text{-non-null}}$ 
with $(\tau^{N} \circ \rho)(A) = 0$ for all $N$, are denoted \emph{$\pi$-reducible} 
(Figure \ref{fig:reducible_set}).  Reducible Markov transitions often partition the ambient 
space into non-overlapping sets in which realized Markov chains are confined 
(Figure \ref{fig:reducible_partition}).  In other words exploration of the ambient space
is ``reduced'' to the local exploration of these divided neighborhoods.

If the $N$-step distribution from every point initialization $x_{0} \in X$, 
\begin{align*}
\tau^{N} \circ \delta_{x_{0}},
\end{align*}
allocates non-zero probability to every $\pi$-non-null set $A$,
\begin{align*}
\tau^{N} \circ \delta_{x_{0}}(A) > 0,
\end{align*}
for some $0 < N(x) < \infty$, then at least some realized Markov chains from every 
initialization will reach every $\pi$-non-null set.  In this case we say that the 
Markov transition is $\pi$-irreducible.

One convenient circumstance where $\pi$-irreducibility is guaranteed is when the 
density of an $N$-step distribution with respect to the target distribution is 
everywhere positive,
\begin{equation*}
\frac{ \mathrm{d} (\tau^{N} \circ \delta_{x_{0}}) }{ \mathrm{d} \pi } (x) > 0
\end{equation*}
for all $x_{0}, x \in X$.

Critically irreducibility is determined not by the precise probability $\pi(A)$ 
allocated to measurable sets but rather the classification of measurable sets into 
null and non-null sets.  If $\pi(A) = 0$ whenever $\phi(A) = 0$ for some 
$\sigma$-finite measure $\phi$ then $\pi$ and $\phi$ will share the same null sets, 
and $\phi$-irreducibility will immediately imply $\pi$-irreducibility.  In this 
case we say that $\phi$ \emph{dominates} $\pi$ or, equivalently, that $\pi$ is
\emph{absolutely continuous} with respect to $\phi$.

Consequently $\pi$-irreducibility is also guaranteed whenever the density of an 
$N$-step distribution with respect to a dominating base measure is everywhere 
positive,
\begin{equation*}
\frac{ \mathrm{d} (\tau^{N} \circ \delta_{x_{0}}) }{ \mathrm{d} \phi } (x) > 0
\end{equation*}
for all $x_{0}, x \in X$.

The same result also holds if only a component of the $N$-step distribution is
absolutely continuous with respect to $\phi$.  More formally if we can write $N$-step 
distribution as a mixture
\begin{equation*}
\tau^{N} \circ \delta_{x_{0}} = \lambda \, \rho_{1} + (1 - \lambda) \, \rho_{2}
\end{equation*}
with $0 < \lambda(x, x_{0}) \le 1$ then
\begin{equation*}
\frac{ \mathrm{d} \rho_{1} }{ \mathrm{d} \phi } (x) > 0 \text{ for all } x_{0}, x \in X,
\end{equation*}
implies that the Markov transition $\tau$ is $\pi$-irreducible even if $\rho_{2}$ 
is not dominated by the base measure $\phi$.  This form is especially convenient 
when working with Metropolis-Hastings methods that mix a continuous proposal transition 
with a singular rejection transition.

\subsection{Periodicity}

Irreducibility ensures that at least some realized Markov chains from any initialization 
will reach every $\pi$-non-null set at least once, but it doesn't guarantee that those sets will 
be visited repeatedly.  Even if a Markov transition is $\pi$-irreducible the exploration 
of some Markov chains can be obstructed by cyclic behavior that forever traps them within 
some subset of the ambient space.  

\begin{figure}
\centering
\begin{tikzpicture}[scale=0.75, thick]
  \draw[white] (-4.25, -4.25) rectangle (4.25, 4.25);
 
  \fill[dark] (0, 1.5) circle (1.25);
  \fill[dark] (-1.5, -1.5) circle (1.25);
  \fill[dark] (1.5, -1.5) circle (1.25);
  \node[dark] at (-1.5, 0.5) { $A_{3}$ };
  \node[dark] at (-3, -2.5) { $A_{1}$ };
  \node[dark] at (0, -2.5) { $A_{2}$ };

  \draw[-{Stealth[length=4pt, width=8pt]}, dark, line width=2] ({0.5 * cos(100)}, {0.5 *  sin(100) - 0.375}) arc (-260:80:0.5);
 
  \pgfmathsetmacro{\N}{29}
 
  \foreach \x/\y [count=\n] in {-3.75/3, -2/2, -1/3.75, 1/3.5, 2/2, 0.453/1.388, -1.635/-1.585, 2.155/-0.666, -0.099/1.812, -1.412/-0.639, 1.515/-1.505, 0.500/0.682} {
    \ifnum \n>1
      \fill[color=mid] (A) circle (4pt);
      \fill[color=light] (A) circle (3pt);
      \draw[color=mid, -{Stealth[length=6pt, width=8pt]}, 
            line width=3, shorten <=2.25pt, shorten >=2pt] (A) -- (\x, \y);
      \draw[color=light, -{Stealth[length=3pt, width=4pt]}, 
            line width=1.5, shorten >=3.25pt] (A) -- (\x, \y);
    \fi
    
    \ifnum \n=\N
      \fill[color=mid] (\x, \y) circle (4pt);
      \fill[color=light] (\x, \y) circle (3pt);
    \fi
    \coordinate (A) at (\x, \y);
  }
    
\end{tikzpicture}

\caption{
A Markov transition is periodic whenever there is a sequence of disjoint, $\pi$-non-null 
sets that trap Markov chains into cyclic transitions.  Here all Markov transitions that 
start in $A_{1}$ are confined to $A_{2}$, those that start in $A_{2}$ are confined to
$A_{3}$, and those that start in $A_{3}$ are confined to $A_{1}$.  Once a Markov chain 
wanders into any of these sets it will be forever doomed to cycle between the three sets 
and unable to explore the rest of the ambient space. 
}
\label{fig:periodicity} 
\end{figure}
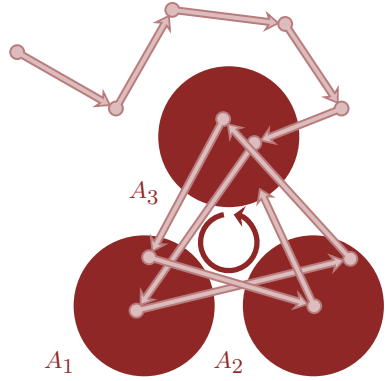

Formally consider a collection of at least two $\pi$-non-null sets,
\begin{equation*} 
A_{1}, \ldots, A_{j}, \ldots, A_{J} \in \mathcal{X}
\end{equation*} 
that are disjoint, 
\begin{equation*}
A_{j} \cap A_{j'} = \emptyset \text{ for } j \ne j'.
\end{equation*}
If $\tau(A_{j + 1} \mid x) = 1$ for all $x \in A_{j}$ then all transitions from points in 
$A_{j}$ will be confined to $A_{j + 1}$.  Moreover if $\tau(A_{1} \mid x) = 1$ for all 
$x \in A_{J}$ then those transitions will eventually return to $A_{1}$ where the cycle 
begins anew (Figure \ref{fig:periodicity}).  In this case the Markov chain is said to 
be \emph{$\pi$-periodic} with \emph{period} $J$ and \emph{periodic decomposition} 
$A_{1}, \ldots, A_{J}$.

A period decomposition serves as a sink that absorbs any Markov chains that venture into 
any of the component sets and preventing them from exploring anywhere else.  In order to 
ensure that Markov chains have the opportunity to explore the entire target distribution 
over and over again we have to avoid this absorbing behavior by ensuring that our Markov 
transition is \emph{$\pi$-aperiodic}.  

One common way to avoid periodic behavior is to ensure that the Markov transition always 
admits a nonzero probability of staying at the initial point.  If $\tau(x \mid x) > 0$ 
then $\tau(A_{j + 1} \mid x) = 1$ only if $x \in A_{j + 1}$.  Consequently transitions
initialized in a set $A_{j}$ cannot be confined to a \emph{disjoint} set $A_{j + 1}$,
obstructing a periodic decomposition.

\subsection{Recurrence}

If a Markov transition with an invariant distribution $\pi$ is both $\pi$-irreducible 
and $\pi$-aperiodic then $\pi$ is the \emph{unique} invariant distribution.  In this case 
the Markov transition is said to be \emph{$\pi$-recurrent}.  The Markov chains generated 
from $\pi$-recurrent Markov transitions will explore every $\pi$-non-null set well enough to 
ensuring asymptotic convergence in the total variation distance,
\begin{equation*}
\lim_{N \rightarrow \infty} | \tau^{N} \circ \delta_{x} - \pi | = 0
\end{equation*}
for $\pi$-almost-all point initializations $x \in X$.  This also implies that the 
empirical average of any measurable function $f : X \rightarrow \mathbb{R}$ with a 
well-defined expectation value $\mathbb{E}_{\pi}[f] < \infty$,
\begin{equation*}
\hat{f}_{N} (x_{0}, \ldots, x_{N})
= \frac{1}{N + 1} \sum_{n = 0}^{N} f \circ \varpi_{n} (x_{0}, \ldots, x_{N})
= \frac{1}{N + 1} \sum_{n = 0}^{N} f(x_{n}) 
\end{equation*}
converges to that exact expectation value,
\begin{equation*}
\lim_{N \rightarrow \infty} \hat{f}_{N} = \mathbb{E}_{\pi}[f]
\end{equation*}
for almost all Markov chain realizations $\{x_{1}, \ldots, x_{n}, \ldots \}$.

Extending this result to all point initializations requires a stronger condition then
$\pi$-recurrence alone.  If for all $\pi$-non-null sets $A \in \mathcal{X}_{\text{non-null}}$ 
and point initializations $x \in X$ we have 
\begin{equation*}
\sum_{n = 1}^{N(x)} (\tau^{N} \circ \delta_{x})(A) = 1
\end{equation*}
for some finite $N(x) < \infty$ then the Markov chain is said to be \emph{Harris recurrent} 
\cite{Harris:1956, Tierney:1994, ChanEtAl:1994}.  Harris recurrence ensures that \emph{all} 
Markov chain realizations from \emph{all} initializations will visit every $\pi$-non-null set 
infinitely often.  This then implies that 
\begin{equation*}
\lim_{N \rightarrow \infty} | \tau^{N} \circ \delta_{x} - \pi | = 0
\end{equation*}
for \emph{all} $x \in X$. 

While a non-zero density of the $N$-step distribution relative to the target distribution,
is a useful sufficient condition for $\pi$-irreducibility, a non-zero density of the $1$-step 
distribution,
\begin{equation*}
\frac{ \mathrm{d} (\tau \circ \delta_{x_{0}}) }{ \mathrm{d} \pi } (x) > 0
\end{equation*}
for all $x_{0}, x \in X$, is sufficient for Harris recurrence.  For any dominating measure 
$\phi$ Harris recurrence is also implied by
\begin{equation*}
\frac{ \mathrm{d} (\tau \circ \delta_{x_{0}}) }{ \mathrm{d} \phi } (x) > 0
\end{equation*}
for all $x_{0}, x \in X$.  As with $\pi$-irreducibility we can also verify Harris 
recurrence using only mixture components of $1$-step distribution to avoid any singular 
behavior.

\section{Preasymptotic Convergence}

Given the inherent difficulty of probabilistic computation we can't take asymptotic convergence 
of Markov chain Monte Carlo estimators for granted.  To avoid problematic behavior we have to 
carefully verify the recurrence properties of a given Markov transition, or rely on Markov
transitions that have been vetted by experts. Unfortunately asymptotic convergence alone does 
not guarantee that Markov chain Monte Carlo will behave well in practice, where we can only ever 
realize \emph{finite} Markov chains.  

Practical Markov chain Monte Carlo performance is instead determined by the \emph{preasymptotic} 
convergence where $N$ is large but finite, which is a much more delicate property.  In this section 
we will review common strategies for quantifying the preasymptotic convergence of Markov chain 
Monte Carlo estimators.

\subsection{Bounding Convergence}

While asymptotic convergence determines \emph{to where} a Markov $N$-step distribution convergences,
preasymptotic convergence determines \emph{how} it gets there.  Typically we quantify preasymptotic 
convergence by bounding the distance between the $N$-step distribution and the target distribution 
as a function of the number of iterations,
\begin{equation*}
|| \tau^{N} \circ \rho - \pi ||_{\mathcal{F}} \le b(\rho, N)
\end{equation*}
for some monotonically decreasing function $b(\rho, N)$ that converges to zero,
\begin{equation*}
\lim_{N \rightarrow \infty} \le b(\rho, N) = 0,
\end{equation*}
If we can bound the distance from every point initialization,
\begin{equation*}
|| \tau^{N} \circ \delta_{x} - \pi ||_{\mathcal{F}} \le b(x, N),
\end{equation*}
then we can bound the distance from any distributional initialization $\rho$ with an expectation value,
\begin{equation*}
b(\rho, N) = \mathbb{E}_{\rho} [ b(-, N) ] \le \underset{x \in X}{\mathrm{sup}} \, b(x, N).
\end{equation*}
Consequently we will focus on bounding the distance between the $N$-step distribution and the target 
distribution from point initializations.

These $N$-step distribution bounds allow us to immediately quantify the preasymptotic behavior of Markov 
chain Monte Carlo estimators.  For example we can work out the preasymptotic bias of the Markov chain 
Monte Carlo estimator for any test function $f \in \mathcal{F}$,
\begin{align*}
\left|
\mathbb{E}_{\tau^{N} \times \delta_{x}} [ \hat{f}_{N} ]
- \mathbb{E}_{\pi}[f]
\right|
&=
\left|
\Big( \frac{1}{N + 1} \sum_{n = 0}^{N} \mathbb{E}_{\tau^{N} \times \delta_{x}} [f \circ \varpi_{n}] \Big)
- \mathbb{E}_{\pi}[f]
\right|
\\
&=
\left|
\frac{1}{N + 1} \sum_{n = 0}^{N} \mathbb{E}_{\tau^{N} \times \delta_{x}}[f \circ \varpi_{n}]
- \frac{1}{N + 1} \sum_{n = 0}^{N} \mathbb{E}_{\pi}[f]
\right|
\\
&=
\left|
\frac{1}{N + 1} \sum_{n = 0}^{N} \mathbb{E}_{\tau^{N} \circ \delta_{x}}[f]
- \frac{1}{N + 1} \sum_{n = 0}^{N} \mathbb{E}_{\pi}[f]
\right|
\\
&=
\frac{1}{N + 1} \left|
\sum_{n = 0}^{N} \Big(
\mathbb{E}_{\tau^{N} \circ \delta_{x}}[f] - \mathbb{E}_{\pi}[f] \Big)
\right|
\\
&\le
\frac{1}{N + 1} \sum_{n = 0}^{N} \Big|
\mathbb{E}_{\tau^{N} \circ \delta_{x}}[f] - \mathbb{E}_{\pi}[f]
\Big|
\\
&\le
\frac{1}{N + 1} \sum_{n = 0}^{N} 
|| \tau^{N} \circ \delta_{x} - \pi ||_{\mathcal{F}}
\\
&\le
\frac{1}{N + 1} \sum_{n = 0}^{N} b(x, N).
\end{align*}
In other words the bias of Markov chain Monte Carlo estimators decays linearly with $N$ 
provided that $\sum_{n = 0}^{N} b(x, N) < N + 1$.  Bounds on the estimator variance 
requires understanding the correlations between Markov chain states and a more delicate 
calculation; see for example \cite{JoulinEtAl:2010}.

If the distance between the $N$-step distribution and target distribution is bounded by a 
geometrically decreasing function of $N$ that is independent of the initialization,
\begin{equation*}
|| \tau^{N} \circ \delta_{x} - \pi ||_{\mathcal{F}} \le b \cdot r^{N}
\end{equation*}
for $0 \le r < 1$, then we say that the Markov transition is \emph{uniformly geometrically ergodic}
in the associated metric.  In this case
\begin{align*}
\left|
\mathbb{E} [ \hat{f}_{N} ]
- \mathbb{E}_{\pi}[f]
\right|
\\
&\le
\frac{1}{N + 1} \sum_{n = 0}^{N} \Big|
\mathbb{E}_{\tau^{n} \circ \rho}[f] - \mathbb{E}_{\pi}[f]
\Big|
\\
&\le
\frac{1}{N + 1} \sum_{n = 0}^{N} b \, r^{N}
\\
&\le
\frac{ b }{N + 1} \sum_{n = 0}^{N} \, r^{N}
\\
&\le
\frac{ b }{N + 1} \frac{ 1 - r^{N + 1} }{1 - r}.
\end{align*}
Because $r < 1$, the bias will monotonically decay at a geometric rate.

Similarly if the metric is bounded by a geometrically decreasing function of $N$
that depends on the initialization,
\begin{equation*}
|| \tau^{N} \circ \rho - \pi ||_{\mathcal{F}} \le b(\rho) r^{N}
\end{equation*}
then we say that the Markov transition is \emph{geometrically ergodic} in the
associated metric.  Here the bias will still decay monotonically, but the distance 
will be affected by the initial condition,
\begin{align*}
\left|
\mathbb{E} [ \hat{f}_{N} ]
- \mathbb{E}_{\pi}[f]
\right|
&\le
\frac{ b(\rho) }{N + 1} \frac{ 1 - r^{N + 1} }{1 - r}.
\end{align*}

Uniform geometric ergodicity is rare and for general problems geometric ergodicity 
is usually the best we can hope to establish.  That said establishing geometric ergodicity 
is no easy task itself, with many bounds readily obstructed by not-uncommon interactions 
between a Markov transition and its stationary distribution.

\subsection{Coupling Methods}

In Section \ref{sec:wasserstein} we used couplings between measures to construct
the $1$-Wasserstein metric.  Couplings between two Markov chains are also useful for 
studying the ergodicity properties of a given Markov transition.

Consider a Markov transition $\tau$ and two different initializations, $\rho$ and
$\omega$, from which we can construct two $N$-step distributions, $\tau^{N} \circ \rho$ 
and $\tau^{N} \circ \omega$.  A \emph{Markov coupling} of $\tau$ is any Markov transition 
$\tau_{\gamma}$ defined on the product space $X \times X$ and initialized from the 
product distribution $\rho \times \omega$ such that every $N$-step distribution 
$\tau_{\gamma}^{N} \circ (\rho \times \omega)$ is a coupling between $\tau^{N} \circ \rho$ 
and $\tau^{N} \circ \omega$.  

A Markov coupling is said to be \emph{contractive} if the distance between the two 
marginal distributions is bounded by a monotonically decreasing function of $N$,
\begin{equation*}
|| \tau^{N} \circ \rho - \tau^{N} \circ \omega ||_{\mathcal{F}} \le b(\rho, \omega, N).
\end{equation*}
While the total variation distance between $N$-step distributions from different 
initializations is always non-expanding the distance from other metrics will be
non-expanding, let alone explicitly contractive, only in special cases.

With clever choices of $\rho$ and $\omega$ a contractive Markov coupling will bound 
the convergence of the Markov transition $\tau$.  For example if we take $\omega = \pi$ 
then
\begin{equation*}
\tau^{N} \circ \omega = \tau^{N} \circ \pi = \pi
\end{equation*}
and
\begin{equation*}
|| \tau^{N} \circ \rho - \pi ||_{\mathcal{F}}
=
|| \tau^{N} \circ \rho - \tau^{N} \circ \omega ||_{\mathcal{F}} 
\le 
b(\rho, \pi, N).
\end{equation*}
In other words a coupling between Markov chain initialized from an arbitrary 
initialization $\rho$ and a Markov chain given a \emph{warm} initializalization
from the stationary distribution can be used to bound the convergence of the 
$N$-step distribution $\tau^{N} \circ \rho$.

\subsubsection{Wasserstein Bounds}

The most common approach to constructing general $1$-Wasserstein bounds is to
consider the $1$-step distributions from two point initializations, $\tau \circ \delta_{x_{1}}$
and $\tau \circ \delta_{x_{2}}$.  Scaling the $1$-Wasserstein distance between these
two distributions by the defining distance function gives the \emph{coarse Ricci curvature}
\cite{Ollivier:2009, JoulinEtAl:2010}
\begin{equation*}
\kappa(x, x') = 1 - \frac{ W_{1, g}(\tau \circ \delta_{x}, \tau \circ \delta_{x'})}{ g(x, x') }
\end{equation*}
If the course Ricci curvature is uniformly lower bounded by some constant,
\begin{equation*}
0 < \kappa \le \kappa(x, x'), \, \forall x, x' \in X
\end{equation*}
then the $1$-Wasserstein distance between the $N$-step distribution from any 
point initialization and the target distribution is bounded by
\begin{align*}
W_{1, g}(\tau^{N} \circ \delta_{x}, \pi)
\le (1 - \kappa)^{N} \cdot  W_{1, g}( \tau \circ \delta_{x}, \pi)
\le (1 - \kappa)^{N} \cdot \mathbb{E}_{\pi}[ g(x,-) ].
\end{align*}

In order to bound the Ricci curvature we have to bound the $1$-Wasserstein distance 
between the $1$-step distributions,
\begin{equation*}
W_{1, g}(\tau \circ \delta_{x}, \tau \circ \delta_{x'}) 
= 
\underset{ \gamma \in \Gamma( \tau \circ \delta_{x}, \tau \circ \delta_{x'} ) }{ \mathrm{inf} }
\mathbb{E}_{\gamma}[g].
\end{equation*}
This, in turn, is bounded by the expectation of the distance function with respect
to \emph{any} coupling between the $1$-step distributions from two different point 
initializations.  Engineering any mathematically-convenient coupling 
$\gamma \in \Gamma(\tau \circ \delta_{x}, \tau \circ \delta_{x'})$ that satisfies
\begin{equation*}
\mathbb{E}_{\gamma}[g] \le g(x, x')
\end{equation*}
for all $x, x' \in X$ immediately establishes an upper bound on the $1$-Wasserstein 
distance,
\begin{equation*}
W_{1, g}(\tau \circ \delta_{x}, \tau \circ \delta_{x'}) \le \mathbb{E}_{\gamma}[g] \le g(x, x'),
\end{equation*}
which then establishes a uniform lower bound on the coarse Ricci curvature,
\begin{equation*}
1 - \kappa = \frac{ W_{1, g}(\tau \circ \delta_{x}, \tau \circ \delta_{x'}) }{ g(x, x') } \le 1.
\end{equation*}
The more strongly Markov chains initialized from two points contract \emph{towards} 
each other in expectation the larger the lower bound we can establish on the Ricci 
curvature and the smaller the upper bound we can establish on the convergence of the 
$N$-step distribution in the $1$-Wasserstein metric.

This expected point-wise contraction is particularly convenient to study when a
Markov transition can be written as a distribution over deterministic transformations 
\cite{DiaconisEtAl:1999}, especially deterministic trajectories.  In this case we
can sometimes demonstrate expected point-wise contraction by coupling the trajectories 
in a way that they tend to evolve towards each other.

\begin{figure}
\centering
\begin{tikzpicture}[scale=0.25]
  \draw[white] (-10, -10) rectangle (10, 10);
 
  \begin{scope}
    \clip (-10, -10) rectangle (10, 10);
    \foreach \A in {-20, -18, ..., 20} {
      \draw[domain={-10:10}, smooth, samples=20, line width=0.5, variable=\x, color=gray70] 
        plot ({\x},{0.005 * \x * \x * \x + 0.1 * \x + \A * cos(10 * \x)});
    }
  \end{scope}
  
  \draw[dark] (0, 0) circle (5 and 7);
  \node[dark] at (3.5, -6.5) { $A$ };

  \pgfmathsetmacro{\xi}{0}
  \pgfmathsetmacro{\xf}{3}
  \pgfmathsetmacro{\A}{4}
  
  \pgfmathsetmacro{\yi}{0.005 * \xi * \xi * \xi + 0.1 * \xi + \A * cos(10 * \xi)}
  \pgfmathsetmacro{\yf}{0.005 * \xf * \xf * \xf + 0.1 * \xf + \A * cos(10 * \xf)}

  \fill[color=mid] (\xi, \yi) circle (0.27);
  
  \draw[domain={\xi:\xf}, smooth, samples=20, variable=\x,
        color=mid, -{Stealth[length=6pt, width=8pt]}, line width=2, shorten >=-1pt] 
        plot ({\x},{0.005 * \x * \x * \x + 0.1 * \x + \A * cos(10 * \x)});
 
  \fill[color=light] (\xi, \yi) circle (0.2);
 
  \draw[domain={\xi:\xf - 0.2}, smooth, samples=20, variable=\x,
        color=light, -{Stealth[length=3pt, width=4pt]}, line width=1, shorten >=-1pt] 
        plot ({\x},{0.005 * \x * \x * \x + 0.1 * \x + \A * cos(10 * \x)});

  \pgfmathsetmacro{\xi}{0}
  \pgfmathsetmacro{\xf}{3}
  \pgfmathsetmacro{\A}{-2}
  
  \pgfmathsetmacro{\yi}{0.005 * \xi * \xi * \xi + 0.1 * \xi + \A * cos(10 * \xi)}
  \pgfmathsetmacro{\yf}{0.005 * \xf * \xf * \xf + 0.1 * \xf + \A * cos(10 * \xf)}

  \fill[color=mid] (\xi, \yi) circle (0.27);
  
  \draw[domain={\xi:\xf}, smooth, samples=20, variable=\x,
        color=mid, -{Stealth[length=6pt, width=8pt]}, line width=2, shorten >=-1pt] 
        plot ({\x},{0.005 * \x * \x * \x + 0.1 * \x + \A * cos(10 * \x)});
 
  \fill[color=light] (\xi, \yi) circle (0.2);
 
  \draw[domain={\xi:\xf - 0.18}, smooth, samples=20, variable=\x,
        color=light, -{Stealth[length=3pt, width=4pt]}, line width=1, shorten >=-1pt] 
        plot ({\x},{0.005 * \x * \x * \x + 0.1 * \x + \A * cos(10 * \x)});


  \pgfmathsetmacro{\xi}{1}
  \pgfmathsetmacro{\xf}{-3}
  \pgfmathsetmacro{\A}{2}
  
  \pgfmathsetmacro{\yi}{0.005 * \xi * \xi * \xi + 0.1 * \xi + \A * cos(10 * \xi)}
  \pgfmathsetmacro{\yf}{0.005 * \xf * \xf * \xf + 0.1 * \xf + \A * cos(10 * \xf)}

  \fill[color=mid] (\xi, \yi) circle (0.27);
  
  \draw[domain={\xi:\xf}, smooth, samples=20, variable=\x,
        color=mid, -{Stealth[length=6pt, width=8pt]}, line width=2, shorten >=-1pt] 
        plot ({\x},{0.005 * \x * \x * \x + 0.1 * \x + \A * cos(10 * \x)});
 
  \fill[color=light] (\xi, \yi) circle (0.2);
 
  \draw[domain={\xi:\xf + 0.18}, smooth, samples=20, variable=\x,
        color=light, -{Stealth[length=3pt, width=4pt]}, line width=1, shorten >=-1pt] 
        plot ({\x},{0.005 * \x * \x * \x + 0.1 * \x + \A * cos(10 * \x)});

  \pgfmathsetmacro{\xi}{1}
  \pgfmathsetmacro{\xf}{-3}
  \pgfmathsetmacro{\A}{-4}
  
  \pgfmathsetmacro{\yi}{0.005 * \xi * \xi * \xi + 0.1 * \xi + \A * cos(10 * \xi)}
  \pgfmathsetmacro{\yf}{0.005 * \xf * \xf * \xf + 0.1 * \xf + \A * cos(10 * \xf)}

  \fill[color=mid] (\xi, \yi) circle (0.27);
  
  \draw[domain={\xi:\xf}, smooth, samples=20, variable=\x,
        color=mid, -{Stealth[length=6pt, width=8pt]}, line width=2, shorten >=-1pt] 
        plot ({\x},{0.005 * \x * \x * \x + 0.1 * \x + \A * cos(10 * \x)});
 
  \fill[color=light] (\xi, \yi) circle (0.2);
 
  \draw[domain={\xi:\xf + 0.18}, smooth, samples=20, variable=\x,
        color=light, -{Stealth[length=3pt, width=4pt]}, line width=1, shorten >=-1pt] 
        plot ({\x},{0.005 * \x * \x * \x + 0.1 * \x + \A * cos(10 * \x)});
 
 
\end{tikzpicture}
\caption{
Within a neighborhood $A$ where the Riemannian sectional curvatures are everywhere 
positive geodesics initialized at nearby points will contract towards each other, 
with the terminal points of the curves closer together than initial points.
}
\label{fig:positive_curvature} 
\end{figure}
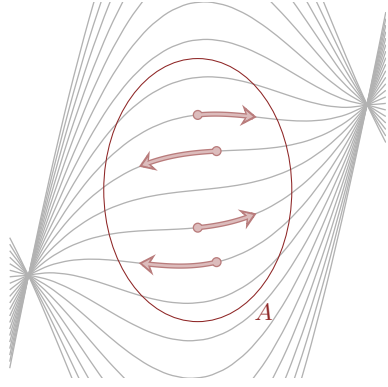

A common example of this approach considers the circumstance where the ambient space 
$X$ is a smooth Riemannian manifold and the Markov transition is a distribution over 
deterministic, geodesic flows.  If the Riemannian sectional curvatures are strictly 
positive in some neighborhood then nearby geodesics will converge towards each other, 
at least for sufficiently short integration times \cite{Lee:2018} 
(Figure \ref{fig:positive_curvature}).  Consequently when the transitions utilize 
sufficiently short geodesics with high enough probability then the local contraction 
of the distance function can imply contraction in expectation, and hence bound the 
$1$-Wasserstein distance of the $1$-step distribution.

Relying on the contraction of deterministic trajectories, however, can be severely 
limiting in practice.  For example in many systems trajectories contract only in 
small neighborhoods, and we will be able to establish the desired bounds only for
Markov transitions that strongly concentrate around the initial point and hence 
Markov chains that explore the target distribution only very slowly.  Alternatively 
we might be able to establish the desired contraction but only for relatively simple 
target distributions, which limits the relevance of the resulting conclusions.

\subsubsection{Total Variation Bounds}

The probability definition of the total variation metric admits a particularly useful 
coupling method for bounding distances between $N$-step distributions and stationary 
distributions \cite{RobertsEtAl:2004}.

A $\pi$-non-null set $A \in \mathcal{X}$ with $\pi(A) > 0$ naturally defines two lifted sets, 
$A \times X \in \mathcal{X} \times \mathcal{X}$ and 
$X \times A \in \mathcal{X} \times \mathcal{X}$.  By definition the probability any 
coupling between any two probability distributions $\mu$ and $\nu$ allocates to these 
sets is just the corresponding marginal probabilities,
\begin{align*}
\gamma(A \times X) &= \mu(A),
\\
\gamma(X \times A) &= \nu(A).
\end{align*}

\begin{figure}
\centering
\begin{tikzpicture}[scale=0.16]

\begin{scope}[shift={(0, 0)}]
  \draw[white] (-15, -14) rectangle(13, 14);
  
  \draw[dark, decorate, decoration={brace,amplitude=5pt}] (-4.5, 10.25) -- (3.5, 10.25);
  \node[dark] at (-0.5, 12.75) { $A$ };
  
  \fill[dark] (-4.5, -10) rectangle (3.5, 10);
  \node[white] at (-0.5, 0) { $A \times X$ };
    
  \draw [->, >=stealth, line width=1] (-10, -10) -- +(20, 0);
  \node[] at (0, -12) { $x_{1}$ };
  \draw [->, >=stealth, line width=1] (-10, -10) -- +(0, 20);
  \node[] at (-12, 0) { $x_{2}$ };
\end{scope}

\begin{scope}[shift={(0, -30)}]
  \draw[white] (-15, -14) rectangle(13, 14);
  
  \draw[dark, decorate, decoration={brace,amplitude=5pt,mirror}] (10.25, -4.5) -- (10.25, 3.5);
  \node[dark] at (12.75, -0.5) { $A$ };
  
  \fill[dark] (-10, -4.5) rectangle (10, 3.5);
  \node[white] at (0, -0.5) { $X \times A$ };
    
  \draw [->, >=stealth, line width=1] (-10, -10) -- +(20, 0);
  \node[] at (0, -12) { $x_{1}$ };
  \draw [->, >=stealth, line width=1] (-10, -10) -- +(0, 20);
  \node[] at (-12, 0) { $x_{2}$ };
\end{scope}

\begin{scope}[shift={(30, 0)}]
  \draw[white] (-15, -14) rectangle(13, 14);
  
  \draw[dark, line width=1.5] (-10, -10) -- (10, 10);
  \node[dark] at (2, -2) { $D$ };
    
  \draw [->, >=stealth, line width=1] (-10, -10) -- +(20, 0);
  \node[] at (0, -12) { $x_{1}$ };
  \draw [->, >=stealth, line width=1] (-10, -10) -- +(0, 20);
  \node[] at (-12, 0) { $x_{2}$ };
\end{scope}

\begin{scope}[shift={(30, -30)}]
  \draw[white] (-15, -14) rectangle(13, 14);
  
  \draw[dark, line width=1.5] (-10, -10) -- (10, 10);
  \node[dark] at (2, -2) { $D$ };
    
  \draw [->, >=stealth, line width=1] (-10, -10) -- +(20, 0);
  \node[] at (0, -12) { $x_{1}$ };
  \draw [->, >=stealth, line width=1] (-10, -10) -- +(0, 20);
  \node[] at (-12, 0) { $x_{2}$ };
\end{scope}

\begin{scope}[shift={(60, 0)}]
  \draw[white] (-15, -14) rectangle(13, 14);
  
  \begin{scope}
    \clip (-10, -4) rectangle (10, 3);
    \draw[dark, line width=1.5] (-10, -10) -- (10, 10);
  \end{scope}
  
  \node[dark] at (3, -6) { $(A \times X) \cap D$ };

  \draw [->, >=stealth, line width=1] (-10, -10) -- +(20, 0);
  \node[] at (0, -12) { $x_{1}$ };
  \draw [->, >=stealth, line width=1] (-10, -10) -- +(0, 20);
  \node[] at (-12, 0) { $x_{2}$ };
\end{scope}

\begin{scope}[shift={(60, -30)}]
  \draw[white] (-15, -14) rectangle(13, 14);
  
  \begin{scope}
    \clip (-10, -4) rectangle (10, 3);
    \draw[dark, line width=1.5] (-10, -10) -- (10, 10);
  \end{scope}
  
  \node[dark] at (3, -6) { $(X \times A) \cap D$ };

  \draw [->, >=stealth, line width=1] (-10, -10) -- +(20, 0);
  \node[] at (0, -12) { $x_{1}$ };
  \draw [->, >=stealth, line width=1] (-10, -10) -- +(0, 20);
  \node[] at (-12, 0) { $x_{2}$ };
\end{scope}

\end{tikzpicture}
\caption{
Any measurable set $A \in \mathcal{X}$ lifts into two measurable sets on 
the product space $X \times X$.  These two lifted sets have exactly the 
same intersection with the diagonal set, $D$.
}
\label{fig:equality_intersections} 
\end{figure}
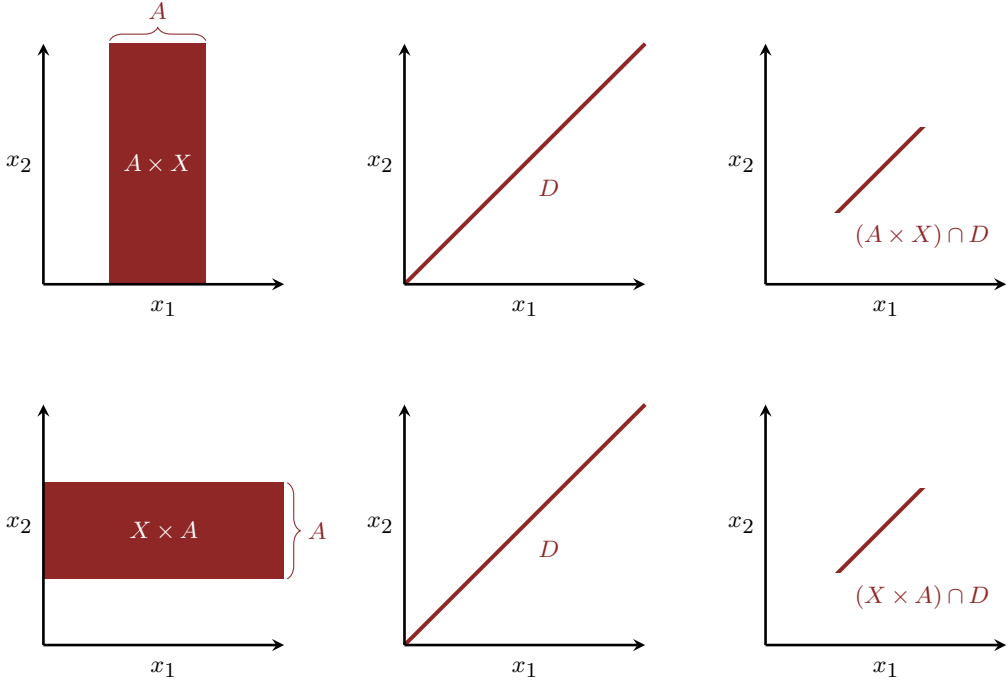

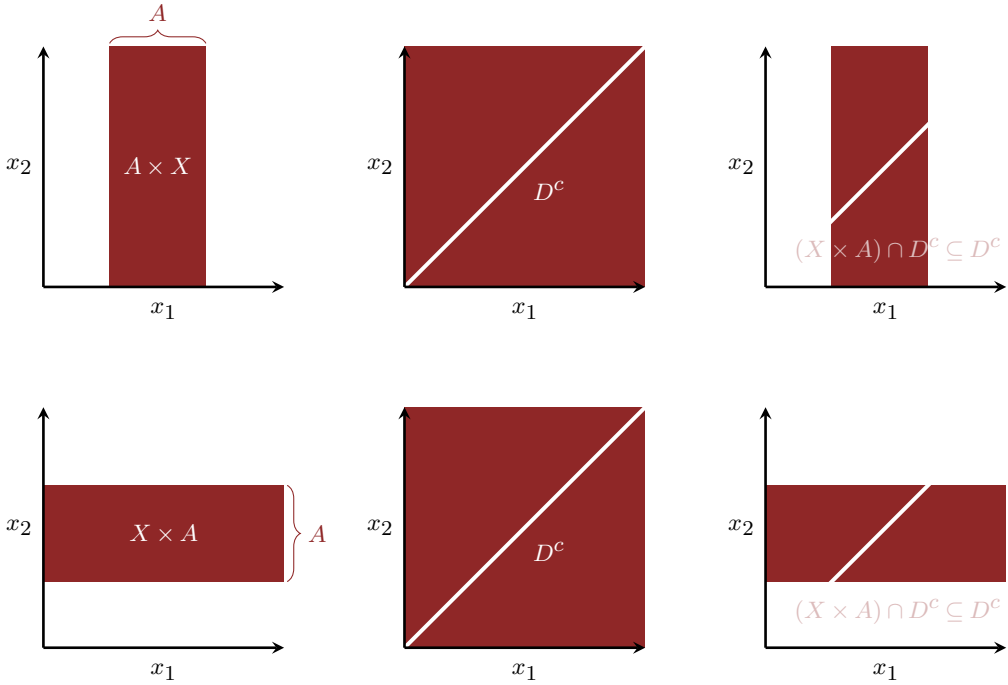
\begin{figure}
\centering
\begin{tikzpicture}[scale=0.16]

\begin{scope}[shift={(0, 0)}]
  \draw[white] (-15, -14) rectangle(13, 14);
  
  \draw[dark, decorate, decoration={brace,amplitude=5pt}] (-4.5, 10.25) -- (3.5, 10.25);
  \node[dark] at (-0.5, 12.75) { $A$ };
  
  \fill[dark] (-4.5, -10) rectangle (3.5, 10);
  \node[white] at (-0.5, 0) { $A \times X$ };
    
  \draw [->, >=stealth, line width=1] (-10, -10) -- +(20, 0);
  \node[] at (0, -12) { $x_{1}$ };
  \draw [->, >=stealth, line width=1] (-10, -10) -- +(0, 20);
  \node[] at (-12, 0) { $x_{2}$ };
\end{scope}

\begin{scope}[shift={(0, -30)}]
  \draw[white] (-15, -14) rectangle(13, 14);
  
  \draw[dark, decorate, decoration={brace,amplitude=5pt,mirror}] (10.25, -4.5) -- (10.25, 3.5);
  \node[dark] at (12.75, -0.5) { $A$ };
  
  \fill[dark] (-10, -4.5) rectangle (10, 3.5);
  \node[white] at (0, -0.5) { $X \times A$ };
    
  \draw [->, >=stealth, line width=1] (-10, -10) -- +(20, 0);
  \node[] at (0, -12) { $x_{1}$ };
  \draw [->, >=stealth, line width=1] (-10, -10) -- +(0, 20);
  \node[] at (-12, 0) { $x_{2}$ };
\end{scope}

\begin{scope}[shift={(30, 0)}]
  \draw[white] (-15, -14) rectangle(13, 14);
  
  \fill[dark] (-10, -10) rectangle (10, 10);
  \draw[white, line width=1.5] (-10, -10) -- (10, 10);
  \node[white] at (2, -2) { $D^{c}$ };
    
  \draw [->, >=stealth, line width=1] (-10, -10) -- +(20, 0);
  \node[] at (0, -12) { $x_{1}$ };
  \draw [->, >=stealth, line width=1] (-10, -10) -- +(0, 20);
  \node[] at (-12, 0) { $x_{2}$ };
\end{scope}

\begin{scope}[shift={(30, -30)}]
  \draw[white] (-15, -14) rectangle(13, 14);
  
  \fill[dark] (-10, -10) rectangle (10, 10);
  \draw[white, line width=1.5] (-10, -10) -- (10, 10);
  \node[white] at (2, -2) { $D^{c}$ };
    
  \draw [->, >=stealth, line width=1] (-10, -10) -- +(20, 0);
  \node[] at (0, -12) { $x_{1}$ };
  \draw [->, >=stealth, line width=1] (-10, -10) -- +(0, 20);
  \node[] at (-12, 0) { $x_{2}$ };
\end{scope}

\begin{scope}[shift={(60, 0)}]
  \draw[white] (-15, -14) rectangle(13, 14);
  
  \fill[dark] (-4.5, -10) rectangle (3.5, 10);
  \draw[white, line width=1.5] (-10, -10) -- (10, 10);
  
  \node[light] at (1, -7) { $(X \times A) \cap D^{c} \subseteq D^{c}$ };
  
  \draw [->, >=stealth, line width=1] (-10, -10) -- +(20, 0);
  \node[] at (0, -12) { $x_{1}$ };
  \draw [->, >=stealth, line width=1] (-10, -10) -- +(0, 20);
  \node[] at (-12, 0) { $x_{2}$ };
\end{scope}

\begin{scope}[shift={(60, -30)}]
  \draw[white] (-15, -14) rectangle(13, 14);
  
  \fill[dark] (-10, -4.5) rectangle (10, 3.5);
  \draw[white, line width=1.5] (-10, -10) -- (10, 10);
  
  \node[light] at (1, -7) { $(X \times A) \cap D^{c} \subseteq D^{c}$ };
  
  \draw [->, >=stealth, line width=1] (-10, -10) -- +(20, 0);
  \node[] at (0, -12) { $x_{1}$ };
  \draw [->, >=stealth, line width=1] (-10, -10) -- +(0, 20);
  \node[] at (-12, 0) { $x_{2}$ };
\end{scope}

\end{tikzpicture}
\caption{
The intersections of the lifted sets $A \times X$ and $X \times A$ with the
inequality set $D^{c}$ are not the same, but they both define subsets of the
inequality set.
}
\label{fig:inequality_intersections} 
\end{figure}

Another natural set to consider when coupling probability distributions together is 
the equality or diagonal set,
\begin{equation*}
D = \{ (x_{1}, x_{2}) \in X \times X \mid x_{1} = x_{2} \},
\end{equation*}
and its complement, the inequality set $D^{c}$.  In particular we can always decompose 
any set on the product space into disjoint intersections with the diagonal set and its 
complement (Figures \ref{fig:equality_intersections}, \ref{fig:inequality_intersections}),
\begin{equation*}
A_{1} \times A_{2} = \big( (A_{1} \times A_{2}) \cap D \big) \cup \big( (A_{1} \times A_{2}) \cap D^{c} \big).
\end{equation*}

This decomposition allows us write the total variation distance between
two probability distributions $\mu$ and $\nu$ as
\begin{align*}
|| \mu - \nu  ||_{\mathrm{TV}} 
&= \underset{ A \in \mathcal{X} }{ \mathrm{sup} } | \mu(A) - \nu(A) |
\\
&= \underset{ A \in \mathcal{X} }{ \mathrm{sup} } | \gamma(A \times X) - \nu(X \times A) |
\\
&= \underset{ A \in \mathcal{X} }{ \mathrm{sup} }
\big| \quad\; \gamma \big( \big( (A \times X) \cap D \big) \cup \big( (A \times X) \cap D^{c} \big) \big)
\\
& \quad\quad\quad\;\; - \gamma \big( \big( (X \times A) \cap D \big) \cup \big( (X \times A) \cap D^{c} \big) \big) \big|
\end{align*}
Because each intersection is disjoint the allocated probabilities decompose into 
independent contributions,
\begin{align*}
|| \mu - \nu  ||_{\mathrm{TV}} 
&= \underset{ A \in \mathcal{X} }{ \mathrm{sup} } 
\big| 
\quad\; \gamma ( (A \times X) \cap D ) + \gamma ( (A \times X) \cap D^{c} )
\\
& \quad\quad\quad\;\; - \gamma ( (X \times A) \cap D ) - \gamma ( (X \times A) \cap D^{c} ) \big|
\\
&= \underset{ A \in \mathcal{X} }{ \mathrm{sup} } 
\big| \quad\; \gamma ( (A \times X) \cap D ) - \gamma ( (X \times A) \cap D )
\\
& \quad\quad\quad\;\;  + \gamma ( (A \times X) \cap D^{c} ) - \gamma ( (X \times A) \cap D^{c} ) \big|.
\end{align*}

The sets $(A \times X) \cap D$ and $(X \times A) \cap D$, however, are exactly the same
and so too must be the probabilities allocated to them,
\begin{align*}
|| \mu - \nu  ||_{\mathrm{TV}} 
&= \underset{ A \in \mathcal{X} }{ \mathrm{sup} } 
\big| \quad\; \gamma ( (A \times X) \cap D ) - \gamma ( (X \times A) \cap D )
\\
& \quad\quad\quad\;\;  + \gamma ( (A \times X) \cap D^{c} ) - \gamma ( (X \times A) \cap D^{c} ) \big|
\\
&= \underset{ A \in \mathcal{X} }{ \mathrm{sup} } 
\big| \hspace{30mm} 0 \hspace{30mm}
\\
& \quad\quad\quad\;\;  + \gamma ( (A \times X) \cap D^{c} ) - \gamma ( (X \times A) \cap D^{c} ) \big|
\\
&= \underset{ A \in \mathcal{X} }{ \mathrm{sup} } 
\big| \gamma ( (A \times X) \cap D^{c} ) - \gamma ( (X \times A) \cap D^{c} ) \big|.
\end{align*}

We can now employ the triangle inequality to separate the absolute value into the
sum of two probabilities,
\begin{align*}
|| \mu - \nu  ||_{\mathrm{TV}} 
&= \underset{ A \in \mathcal{X} }{ \mathrm{sup} } 
| \gamma ( (A \times X) \cap D^{c} ) - \gamma ( (X \times A) \cap D^{c} ) |
\\
&\le \underset{ A \in \mathcal{X} }{ \mathrm{sup} } 
| \gamma ( (A \times X) \cap D^{c} ) | + | \gamma ( (X \times A) \cap D^{c} ) |.
\end{align*}

Because the remaining intersection sets are both subsets of $D^{c}$,
\begin{align*}
(A \times X) \cap D^{c} &\subseteq D^{c}
\\
(X \times A) \cap D^{c} &\subseteq D^{c},
\end{align*}
we have
\begin{align*}
\gamma ( (A \times X) \cap D^{c} ) &\le \gamma(D^{c})
\\
\gamma ( (X \times A) \cap D^{c} ) &\le \gamma(D^{c}),
\end{align*}
and
\begin{align*}
|| \mu - \nu  ||_{\mathrm{TV}} 
&\le \underset{ A \in \mathcal{X} }{ \mathrm{sup} } \,
| \gamma ( (A \times X) \cap D^{c} ) | + | \gamma ( (A \times X) \cap D^{c} ) |
\\
&\le \underset{ A \in \mathcal{X} }{ \mathrm{sup} } \, 2 \, | \gamma(D^{c}) |
\\
&\le \underset{ A \in \mathcal{X} }{ \mathrm{sup} }\, 2 \, \gamma(D^{c})
\end{align*}

Finally because this bound is independent of $A$ we can drop the supremum to 
give
\begin{align*}
|| \mu - \nu  ||_{\mathrm{TV}} 
&\le \underset{ A \in \mathcal{X} }{ \mathrm{sup} } \, \gamma(D^{c})
\\
&\le \gamma(D^{c}).
\end{align*}
Consequently if we can bound the probability allocated to the inequality set for 
\emph{any} coupling between $\mu$ and $\nu$ then we can bound the total variation 
distance between the two probability distributions.

\begin{figure}
\centering
\begin{tikzpicture}[scale=0.75, thick]

\begin{scope}[shift={(0, 0)}]
  \draw[dark, dashed] (-4.25, -4.25) rectangle (4.25, 4.25);
 
  \fill[dark] (0, 0) circle (1);
  \node[dark] at (1.5, 0) { $C$ };
 
  \pgfmathsetmacro{\N}{8}
 
  \foreach \x/\y [count=\n] in {-2.829/1.958, -2.079/3.108, -1.406/2.073, -1.060/0.770, -0.576/-0.163, -0.091/-1.526, 0.758/-2.528, -0.479/-3.533} {
    \ifnum \n>1
      \fill[color=mid] (A) circle (4pt);
      \fill[color=light] (A) circle (3pt);
      \draw[color=mid, -{Stealth[length=6pt, width=8pt]}, 
            line width=3, shorten <=2.25pt, shorten >=2pt] (A) -- (\x, \y);
      \draw[color=light, -{Stealth[length=3pt, width=4pt]}, 
            line width=1.5, shorten >=3.25pt] (A) -- (\x, \y);
    \fi
    
    \ifnum \n=\N
      \fill[color=mid] (\x, \y) circle (4pt);
      \fill[color=light] (\x, \y) circle (3pt);
    \fi
    \coordinate (A) at (\x, \y);
  }
  
  \foreach \x/\y [count=\n] in {3.280/3.561, 1.989/2.797, 1.639/2.027, 1.221/1.401, 0.433/0.391, -0.091/-1.526} {
    \ifnum \n>1
      \fill[color=mid] (A) circle (4pt);
      \fill[color=light] (A) circle (3pt);
      \draw[color=mid, -{Stealth[length=6pt, width=8pt]}, 
            line width=3, shorten <=2.25pt, shorten >=2pt] (A) -- (\x, \y);
      \draw[color=light, -{Stealth[length=3pt, width=4pt]}, 
            line width=1.5, shorten >=3.25pt] (A) -- (\x, \y);
    \fi
    
    \ifnum \n=\N
      \fill[color=mid] (\x, \y) circle (4pt);
      \fill[color=light] (\x, \y) circle (3pt);
    \fi
    \coordinate (A) at (\x, \y);
  }
\end{scope}

\begin{scope}[shift={(10, 0)}]
  \draw[dark, dashed] (-4.25, -4.25) rectangle (4.25, 4.25);
 
  \fill[dark] (0, 0) circle (1);
  \node[dark] at (1.5, 0) { $C$ };
 
  \pgfmathsetmacro{\N}{8}
 
  \foreach \x/\y [count=\n] in {-2.829/1.958, -2.079/3.108, -1.406/2.073, -1.060/0.770, -0.576/-0.163, -0.931/-0.966, -2.139/-2.335, -2.644/-2.102} {
    \ifnum \n>1
      \fill[color=mid] (A) circle (4pt);
      \fill[color=light] (A) circle (3pt);
      \draw[color=mid, -{Stealth[length=6pt, width=8pt]}, 
            line width=3, shorten <=2.25pt, shorten >=2pt] (A) -- (\x, \y);
      \draw[color=light, -{Stealth[length=3pt, width=4pt]}, 
            line width=1.5, shorten >=3.25pt] (A) -- (\x, \y);
    \fi
    
    \ifnum \n=\N
      \fill[color=mid] (\x, \y) circle (4pt);
      \fill[color=light] (\x, \y) circle (3pt);
    \fi
    \coordinate (A) at (\x, \y);
  }
  
  \foreach \x/\y [count=\n] in {3.280/3.561, 1.989/2.797, 1.639/2.027, 1.221/1.401, 0.433/0.391, 1.102/-1.323, 1.871/-2.202, 1.971/-2.969} {
    \ifnum \n>1
      \fill[color=mid] (A) circle (4pt);
      \fill[color=light] (A) circle (3pt);
      \draw[color=mid, -{Stealth[length=6pt, width=8pt]}, 
            line width=3, shorten <=2.25pt, shorten >=2pt] (A) -- (\x, \y);
      \draw[color=light, -{Stealth[length=3pt, width=4pt]}, 
            line width=1.5, shorten >=3.25pt] (A) -- (\x, \y);
    \fi
    
    \ifnum \n=\N
      \fill[color=mid] (\x, \y) circle (4pt);
      \fill[color=light] (\x, \y) circle (3pt);
    \fi
    \coordinate (A) at (\x, \y);
  }
\end{scope}
  
\end{tikzpicture}
\caption{
The splitting Markov coupling evolves two Markov chains independently until they reach a 
small set, $C \in \mathcal{X}$.  Once in the small set the splitting Markov coupling mixes
(left) a transition that merges the two marginal Markov transitions together for all future 
iterations and (right) a transition that continues to evolve them independently.  The
probability of merging bounds the probability allocated to the inequality set, 
$\gamma(D^{c})$.
}
\label{fig:splitting_coupling} 
\end{figure}
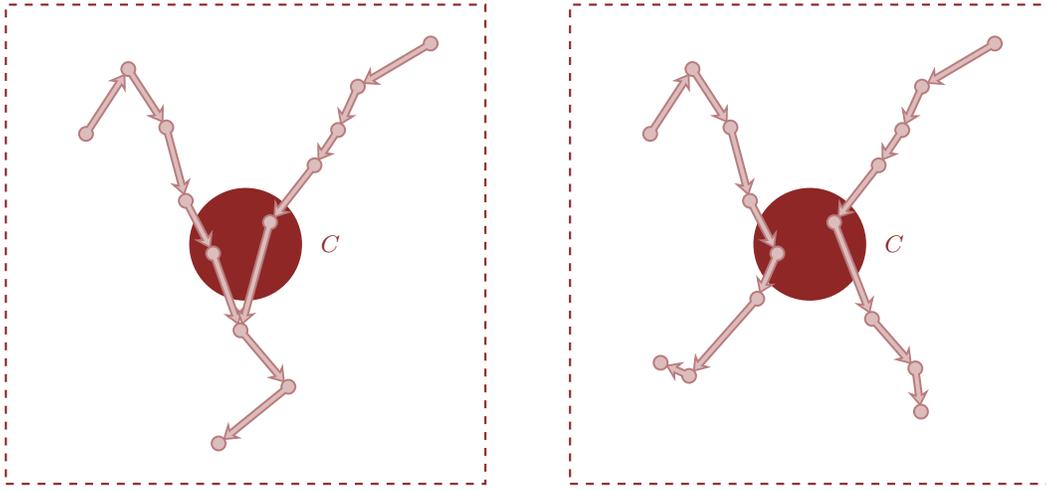

Under certain conditions we can construct a Markov coupling between $\tau^{N} \circ \rho$ and 
$\tau^{N} \circ \pi$ that concentrates enough probability on $D^{c}$ to ensure a finite bound.
Conceptual this \emph{splitting coupling} evolves the marginal Markov chains independently
until they both fall into a \emph{small set} $C \in X$ at the same iteration.  At this 
point the Markov chains merge with a certain probability; if they don't merge then they 
continue to transition independently until the next meeting, but if they do merge then 
they transition to the same point and share the exact same transition for all future 
iterations (Fig \ref{fig:splitting_coupling}).

The probability $\gamma(D^{c})$ allocated by this coupling is determined by probability 
that the marginal Markov transitions $\tau^{N} \circ \rho$ and $\tau^{N} \circ \pi$ merge once 
in the small set.  In order to ensure that this probability is non-zero the Markov transition 
being coupled needs to admit \emph{minorization} and \emph{drift} conditions.

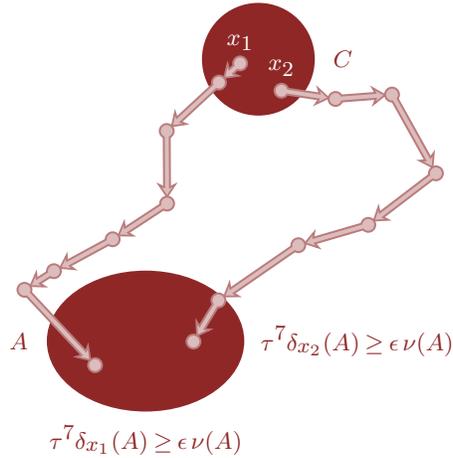
\begin{figure}
\centering
\begin{tikzpicture}[scale=0.75, thick]

  \draw[white] (-4.75, -7.25) rectangle (3.70, 1.25);
 
  \fill[dark] (0, 0) circle (1);
  \node[dark] at (1.5, 0) { $C$ };
  \node[white] at (-0.324, 0.3) { $x_{1}$ };
  \node[white] at (0.409, -0.15) { $x_{2}$ };
 
  \fill[dark] (-2, -5) circle (1.75 and 1.25);
  \node[dark] at (-4.25, -5) { $A$ };
  \node[dark] at (1.75, -5) { $\tau^{7} \delta_{x_{2}}(A) \ge \epsilon \, \nu(A)$ };
  \node[dark] at (-2, -6.75) { $\tau^{7} \delta_{x_{1}}(A) \ge \epsilon \, \nu(A)$ };
 
  \pgfmathsetmacro{\N}{8}
 
  \foreach \x/\y [count=\n] in {-0.324/-0.074, -0.699/-0.414, -1.625/-1.273, -1.618/-2.561, -2.583/-3.194, -3.624/-3.760, -4.146/-4.090, -2.895/-5.427} {
    \ifnum \n>1
      \fill[color=mid] (A) circle (4pt);
      \fill[color=light] (A) circle (3pt);
      \draw[color=mid, -{Stealth[length=6pt, width=8pt]}, 
            line width=3, shorten <=2.25pt, shorten >=2pt] (A) -- (\x, \y);
      \draw[color=light, -{Stealth[length=3pt, width=4pt]}, 
            line width=1.5, shorten >=3.25pt] (A) -- (\x, \y);
    \fi
    
    \ifnum \n=\N
      \fill[color=mid] (\x, \y) circle (4pt);
      \fill[color=light] (\x, \y) circle (3pt);
    \fi
    \coordinate (A) at (\x, \y);
  }
  
  \foreach \x/\y [count=\n] in {0.409/-0.558, 1.367/-0.705, 2.369/-0.627, 3.147/-2.024, 1.948/-2.938, 0.708/-3.299, -0.708/-4.278, -1.150/-5.009} {
    \ifnum \n>1
      \fill[color=mid] (A) circle (4pt);
      \fill[color=light] (A) circle (3pt);
      \draw[color=mid, -{Stealth[length=6pt, width=8pt]}, 
            line width=3, shorten <=2.25pt, shorten >=2pt] (A) -- (\x, \y);
      \draw[color=light, -{Stealth[length=3pt, width=4pt]}, 
            line width=1.5, shorten >=3.25pt] (A) -- (\x, \y);
    \fi
    
    \ifnum \n=\N
      \fill[color=mid] (\x, \y) circle (4pt);
      \fill[color=light] (\x, \y) circle (3pt);
    \fi
    \coordinate (A) at (\x, \y);
  }
  
\end{tikzpicture}
\caption{
A minorization condition ensures that the Markov transitions from any two point initializations
in the small set $C \in \mathcal{X}$ have a non-zero overlap on any $\pi$-non-null set $A$ after a 
finite number of iterations.  In other words once in the small set all Markov transitions look 
somewhat similar, which allows for two independent Markov chains in a splitting coupling to 
merge with non-zero probability.
}
\label{fig:minorization} 
\end{figure}

A Markov transition admits a minorization condition if there exists a positive real number 
$\epsilon > 0$, a positive integer $M$, and measure $\nu$ that dominates $\pi$ such 
that
\begin{equation*}
(\tau^{M} \circ \delta_{x}) (A) \ge \epsilon \, \nu(A)
\end{equation*}
for all $A \in \mathcal{X}$ and points $x$ in the set $C \subset X$ (Figure \ref{fig:minorization}).  
Conceptually the minorization condition ensures that the $M$-step distributions from 
all point initializations in the small set, $\tau^{M} \circ \delta_{x}$, share a finite overlap.
This overlap then allows two coupled Markov transitions to merge with finite probability 
once in the small set.

If a minorization condition with overlap $\epsilon$ holds after $M$ iterations then the 
probability that the splitting coupling does \emph{not} merge two Markov chains initialized 
within the small set after $M$ transitions is no greater than $(1 - \epsilon)^{M}$.  This then 
ensures the bound
\begin{equation*}
|| \tau^{M} \circ \delta_{x_{1}} - \tau^{M} \circ \delta_{x_{2}} ||_{\mathrm{TV}} 
\le \gamma(D^{c}) 
\le (1 - \epsilon)^{M}
\end{equation*}
for $x_{1}, x_{2} \in C$.

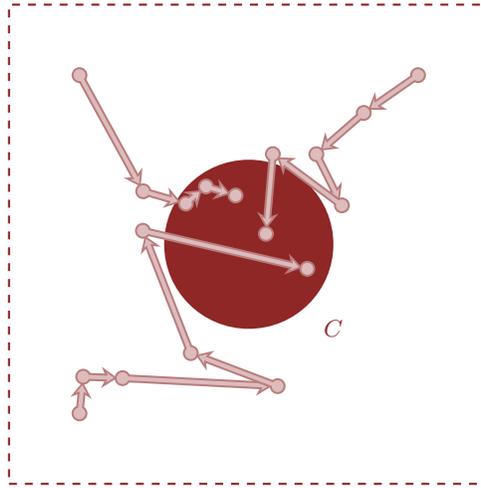
\begin{figure}
\centering
\begin{tikzpicture}[scale=0.75, thick]

\draw[dark, dashed] (-4.25, -4.25) rectangle (4.25, 4.25);
  \fill[dark] (0, 0) circle (1.5);
  \node[dark] at (1.5, -1.5) { $C$ };
 
  \pgfmathsetmacro{\N}{5}
  
  \foreach \x/\y [count=\n] in {-3.000/3.000, -1.871/0.940, -1.114/0.721, -0.758/1.025, -0.234/0.865} {
    \ifnum \n>1
      \fill[color=mid] (A) circle (4pt);
      \fill[color=light] (A) circle (3pt);
      \draw[color=mid, -{Stealth[length=6pt, width=8pt]}, 
            line width=3, shorten <=2.25pt, shorten >=2pt] (A) -- (\x, \y);
      \draw[color=light, -{Stealth[length=3pt, width=4pt]}, 
            line width=1.5, shorten >=3.25pt] (A) -- (\x, \y);
    \fi
    
    \ifnum \n=\N
      \fill[color=mid] (\x, \y) circle (4pt);
      \fill[color=light] (\x, \y) circle (3pt);
    \fi
    \coordinate (A) at (\x, \y);
  }
  
  \pgfmathsetmacro{\N}{6}
  
  \foreach \x/\y [count=\n] in {3.000/3.000, 2.038/2.330, 1.198/1.598, 1.651/0.690, 0.427/1.600, 0.304/0.185} {
    \ifnum \n>1
      \fill[color=mid] (A) circle (4pt);
      \fill[color=light] (A) circle (3pt);
      \draw[color=mid, -{Stealth[length=6pt, width=8pt]}, 
            line width=3, shorten <=2.25pt, shorten >=2pt] (A) -- (\x, \y);
      \draw[color=light, -{Stealth[length=3pt, width=4pt]}, 
            line width=1.5, shorten >=3.25pt] (A) -- (\x, \y);
    \fi
    
    \ifnum \n=\N
      \fill[color=mid] (\x, \y) circle (4pt);
      \fill[color=light] (\x, \y) circle (3pt);
    \fi
    \coordinate (A) at (\x, \y);
  }
  
  \pgfmathsetmacro{\N}{7}
  
  \foreach \x/\y [count=\n] in {-3.000/-3.000, -2.934/-2.350, -2.237/-2.375, 0.512/-2.517, -1.028/-1.931, -1.876/0.239, 1.037/-0.435} {
    \ifnum \n>1
      \fill[color=mid] (A) circle (4pt);
      \fill[color=light] (A) circle (3pt);
      \draw[color=mid, -{Stealth[length=6pt, width=8pt]}, 
            line width=3, shorten <=2.25pt, shorten >=2pt] (A) -- (\x, \y);
      \draw[color=light, -{Stealth[length=3pt, width=4pt]}, 
            line width=1.5, shorten >=3.25pt] (A) -- (\x, \y);
    \fi
    
    \ifnum \n=\N
      \fill[color=mid] (\x, \y) circle (4pt);
      \fill[color=light] (\x, \y) circle (3pt);
    \fi
    \coordinate (A) at (\x, \y);
  }  
\end{tikzpicture}
\caption{
A drift condition ensures that Markov chains steadily converge towards smaller values of 
the associated drift function.  When the drift function is minimized within a small set 
$C$ then Markov chains will steadily converge towards $C$, and multiple Markov chains
have a finite probability of meeting in $C$ at the same time.
}
\label{fig:drift} 
\end{figure}

This minorization bound is only relevant, however, if two Markov chains are able to meet in the 
small set in the first place.  In particular the splitting Markov coupling needs to assign 
sufficient probability to the lifted set $C \times C$.  One way to guarantee a finite meeting 
probability is a \emph{drift condition} which requires the existence of a \emph{drift function} 
$V : X \rightarrow [1, \infty) \subset \mathbb{R}$ satisfying
\begin{equation*}
\mathbb{E}_{\tau \circ \delta_{x}} [ V ] \le \lambda \, V(x) + b \, \mathbb{I}_{C}(x)
\end{equation*}
for constants $0 < \lambda < 1$ and $b < \infty$ and all $x \in X$.  If the the minimum of this 
drift function is within the small set $C$ then a drift condition ensures that Markov chains
initialized outside of the small set will transition towards the small set with high probability
(Figure \ref{fig:drift}).  This then ensures that the probability of two independent Markov 
chains not meeting in the small set after $N$ iterations is bounded by a geometric function of 
$N$.

In other words the drift condition ensures a geometric bound on two Markov chains
not meeting in the small set, while the minorization condition ensures a bound on 
$\gamma(D^{c})$ once the Markov chains have met.  Together with $\phi$-irreducibility
these conditions guarantee a geometric bound on the total variation distance between
the two $N$-step distributions,
\begin{equation*}
|| \tau^{N} \circ \delta_{x} - \tau^{N} \circ \omega ||_{\mathrm{TV}} \le C \, V(x) r^{N},
\end{equation*}
and hence the convergence of any $N$-step distribution towards the target distribution,
\begin{equation*}
|| \tau^{N} \circ \delta_{x} - \pi ||_{\mathrm{TV}}
=
|| \tau^{N} \circ \delta_{x} - \tau^{N} \circ \pi ||_{\mathrm{TV}}
\le C \, V(x) r^{N}.
\end{equation*}

One interesting consequence of this construction is that the drift condition ensures 
not only a geometric convergence bound in the total variation distance but also the 
$V$-norm,
\begin{equation*}
|| \tau^{N} \circ \delta_{x} - \pi ||_{\mathcal{F}_{V}}
\end{equation*}
where
\begin{equation*}
\mathcal{F}_{V} = \{ f \in C^{0}(X) \mid | f(x) | < V(x) \}.
\end{equation*}
The faster the drift function grows the larger this space of test functions will be
and the more applicable it might be in practice.

Unfortunately even when the space of test functions is large the explicit convergence 
bound guaranteed by this splitting coupling construction is often extremely loose. 
The configuration of a Markov chain Monte Carlo algorithm motivated by an explicit 
bound, such as for how many iterations we need to run each Markov chain to achieve 
a given estimator error, is often drastically conservative. 

\subsection{The Markov Chain Monte Carlo Central Limit Theorem}

Although geometric total variation bounds might not be directly useful for quantifying 
the preasymptotic behavior of Markov transitions, they are still extremely useful 
in practice.  In particular geometric ergodicity in the total variation metric 
guarantees the existence of a \emph{Markov chain Monte Carlo central limit theorem} 
that allows us to \emph{empirically} quantify preasymptotic convergence of Markov 
chain Monte Carlo estimators for any sufficiently integrable function 
\cite{RobertsEtAl:2004}.

Let $f : X \rightarrow \mathbb{R}$ be a real-valued, measurable function with 
$\mathbb{E}_{\pi}[f^{\delta}]$ finite for any $\delta > 0$.  For any realized Markov chain 
$(x_{0}, \ldots, x_{N})$, we can construct a corresponding Markov chain Monte Carlo 
estimator
\begin{equation*}
\hat{f}_{N}(x_{0}, \ldots, x_{N}) 
= \frac{1}{N + 1} \sum_{n = 0}^{N} f \circ \varpi_{n} (x_{0}, \ldots, x_{N})
= \frac{1}{N + 1} \sum_{n = 0}^{N} f(x_{n}),
\end{equation*}
If the Markov transition is Harris recurrent then the distribution of estimator values 
derived from the possible Markov chain realizations initialized from $\rho$,
$( \hat{f}_{N} )_{*} (\tau^{N} \times \rho)$, asymptotically converges to a Dirac 
distribution at the exact expectation value,
\begin{equation*}
\lim_{N \rightarrow \infty} ( \hat{f}_{N} )_{*} (\tau^{N} \times \rho) = \delta_{ \mathbb{E}_{\pi}[f] }.
\end{equation*}
Moreover if $\tau$ is geometrically ergodic in the total variation distance then the
probability density function of $( \hat{f}_{N} )_{*} (\tau^{N} \times \rho)$ satisfies 
a central limit theorem,
\begin{equation*}
\lim_{N \rightarrow \infty} 
\frac{ ( \hat{f}_{N} )_{*} (\tau^{N} \times \rho) (y) - \mathbb{E}_{\pi}[f] }{ \sigma_{N} } 
= \text{normal}(y \mid 0, 1),
\end{equation*}
where $\sigma_{N}^{2}$ is the \emph{asymptotic variance}.

To give this asymptotic variance an explicit form let $\mu_{f} = \mathbb{E}_{\pi}[f]$ and 
$\sigma^{2}_{f} = \mathbb{E}_{\pi}[ (f - \mu_{f})^{2} ]$ and define the lag-$l$ 
autocovariance of a Markov transition as
\begin{equation*}
\xi_{l}[f] =
\mathbb{E}_{ \tau^{N} \times \rho } [ ( f \circ \varpi_{l + n} - \mu_{f} ) ( f \circ \varpi_{n} - \mu_{f}) ],
\end{equation*}
with the corresponding lag-$l$ autocorrelation defined as
\begin{equation*}
\zeta_{l}[f] = \frac{ \xi_{l}[f] }{ \mathrm{Var}_{\pi}[f] }.
\end{equation*}
The asymptotic variance is then given by
\begin{align*}
\lim_{N \rightarrow \infty} N \cdot \sigma_{N}^{2}
&=
\sum_{l = -\infty}^{\infty} \xi_{l}[f]
\\
&=
\mathrm{Var}_{\pi}[f] \sum_{l = -\infty}^{\infty} \zeta_{l}[f]
\\
&=
\mathrm{Var}_{\pi}[f] \cdot \left( 1 + 2 \cdot \sum_{l = 1}^{\infty} \zeta_{l}[f] \right).
\end{align*}

Although there is an initial bias in Markov chain Monte Carlo estimators from geometrically 
ergodic Markov transitions, that bias decays linearly with the number of iterations $N$.
Consequently for large enough $N$ the bias becomes negligible and we can use the central 
limit theorem to approximate the probability density function of estimator values as
\begin{equation*}
( \hat{f}_{N} )_{*} (\tau^{N} \times \rho) (y) = \text{normal}(y \mid \mathbb{E}_{\pi}[f], \sigma_{N}).
\end{equation*}
For fixed $N$ this defines a \emph{Markov chain Monte Carlo standard error} that 
quantifies how well a realized value of $\hat{f}_{N}$ estimates the exact expectation 
value.  In practice the error is typically written as
\begin{align*}
N \cdot \sigma_{N}^{2}
&=
\mathrm{Var}_{\pi}[f] \cdot \left( 1 + 2 \cdot \sum_{l = 1}^{\infty} \zeta_{l}[f] \right)
\\
\sigma_{N}
&=
\sqrt{ \mathrm{Var}_{\pi}[f] \cdot \frac{ \left( 1 + 2 \cdot \sum_{l = 1}^{\infty} \zeta_{l}[f] \right) }{ N } }
\\
\sigma_{N}
&=
\sqrt{ \mathrm{Var}_{\pi}[f] \cdot \frac{1}{ \text{ESS}[f] } },
\end{align*}
where the \emph{effective sample size}
\begin{equation*}
\text{ESS}[f] = \frac{ N }{ (1 + 2 \cdot \sum_{l = 1}^{\infty} \zeta_{l}[f] )}
\end{equation*}
moderates the precision of the estimator.  The less the values of $f$ are autocorrelated 
in the realized Markov chains the larger the effective sample size, and the more precise 
Markov chain Monte Carlo estimators, will be for any fixed Markov chain length $N$.

In practice if we don't know $\mathbb{E}_{\pi}[f]$ then it's highly unlikely that we will 
know the exact variance or effective sample size, either.  If $\mathbb{E}_{\pi}[f^{4 + \delta}]$
is finite, however, then we can use a realized Markov chain to estimate these quantities 
and hence the Markov chain standard error.  The error in this estimation is comparable to 
the error introduced by assuming that the central limit theorem holds for finite $N$, 
and consequently negligible when $N$ is large enough.

When we know that a Markov transition is geometrically ergodic in the total variation distance
then the implementation of Markov chain Monte Carlo is straightforward.  Given some initialization
distribution $\rho$ we generate a Markov chain by sampling an initial point from $\rho$ and
then applying the Markov transition $N$ times.  For any function $f : X \rightarrow \mathbb{R}$ 
we can construct a Markov chain Monte Carlo estimator for the corresponding expectation value, as 
well as the variance and the autocorrelations needed to estimate the corresponding error 
estimate.

\begin{figure}
\centering
\begin{tikzpicture}[scale=0.25, thick]
  
\begin{scope}[shift={(0, 0)}]

  \draw[white] (-16, -4) rectangle (12, 22);

  \draw[gray60, dashed, line width=1.5] (0, 0)  -- (0, 14);
  \node[black] at (0, 15) { $\mathbb{E}_{\pi}[f]$ };

  \draw[domain={-10:9.8}, smooth, samples=150, line width=1, variable=\x, color=dark] 
    plot ({\x},{10 * laplace(0.75 * (\x - 2)) + 5 * normal(0.9 * (\x + 0)) + 5 * normal(0.7 * (\x + 3))});

  \draw [->, >=stealth, line width=1] (-10, 0) -- (10, 0);
  \node[] at (0, -2) { $y$ };

  \draw [->, >=stealth, line width=1] (-10, 0) -- (-10, 15);
  \node[rotate=90] at (-12.5, 7.5) { $(\hat{f}_{N})_{*} (\tau^{N} \times \rho) (y)$ };

  \node[] at (0, 20) { $N \approx 0$ };
  \node[] at (0, 18) { Initialization Regime };
  
\end{scope}

\begin{scope}[shift={(30, 0)}]

  \draw[white] (-16, -4) rectangle (12, 22);

  \draw[gray60, dashed, line width=1.5] (0, 0)  -- (0, 14);
  \node[black] at (0, 15) { $\mathbb{E}_{\pi}[f]$ };

  \foreach \n in {1, 2, ..., 5} {
    \pgfmathsetmacro{\prop}{20 * \n};
    \colorlet{custom}{dark!\prop!white};
    
    \pgfmathsetmacro{\sigma}{3 / \n};
    \draw[domain={-10:9.8}, smooth, samples=150, line width=1, variable=\x, color=custom] 
      plot ({\x},{6 * normal(\x / \sigma) / \sigma});
  }
  
  \draw [->, >=stealth, line width=1] (-10, 0) -- (10, 0);
  \node[] at (0, -2) { $y$ };

  \draw [->, >=stealth, line width=1] (-10, 0) -- (-10, 15);
  \node[rotate=90] at (-12.5, 7.5) { $(\hat{f}_{N})_{*} (\tau^{N} \times \rho) (y)$ };

  \node[] at (0, 20) { $N \gg 0$ };
  \node[] at (0, 18) { Central Limit Theorem Regime };
  
\end{scope}

\begin{scope}[shift={(15, -28)}]

  \draw[white] (-16, -4) rectangle (12, 22);

  \draw[gray60, dashed, line width=1.5] (0, 0)  -- (0, 14);
  \node[black] at (0, 15) { $\mathbb{E}_{\pi}[f]$ };

  \foreach \n in {1, 2, ..., 5} {
    \pgfmathsetmacro{\prop}{20 * \n};
    \colorlet{custom}{dark!\prop!white};
    
    \pgfmathsetmacro{\sigma}{0.1 / \n};
    \draw[domain={-1:1}, smooth, samples=150, line width=1, variable=\x, color=custom] 
      plot ({\x},{14.5 * normal(\x / \sigma)});
  }
  
  \draw [->, >=stealth, line width=1] (-10, 0) -- (10, 0);
  \node[] at (0, -2) { $y$ };

  \draw [->, >=stealth, line width=1] (-10, 0) -- (-10, 15);
  \node[rotate=90] at (-12.5, 7.5) { $(\hat{f}_{N})_{*} (\tau^{N} \times \rho) (y)$ };

  \node[] at (0, 20) { $N \rightarrow \infty$ };
  \node[] at (0, 18) { Asymptotic Limit };
  
\end{scope}
  
\end{tikzpicture}
\caption{
Markov chain Monte Carlo estimators from geometrically ergodic Markov transitions 
converge in three distinct phases.  In the initial phase where the number of iterations
small, $N \approx 0$, the estimators are biased towards the initialization.  As the 
number of iterations grows this bias decays, and for large enough iterations $N \gg 0$
the distribution of estimator values is well-approximated by a normal density function
centered on the exact expectation value.  After an infinite number of iterations this
normal approximation collapses to a Dirac distribution centered on the exact expectation 
value.
}
\label{fig:phases_of_convergence} 
\end{figure}
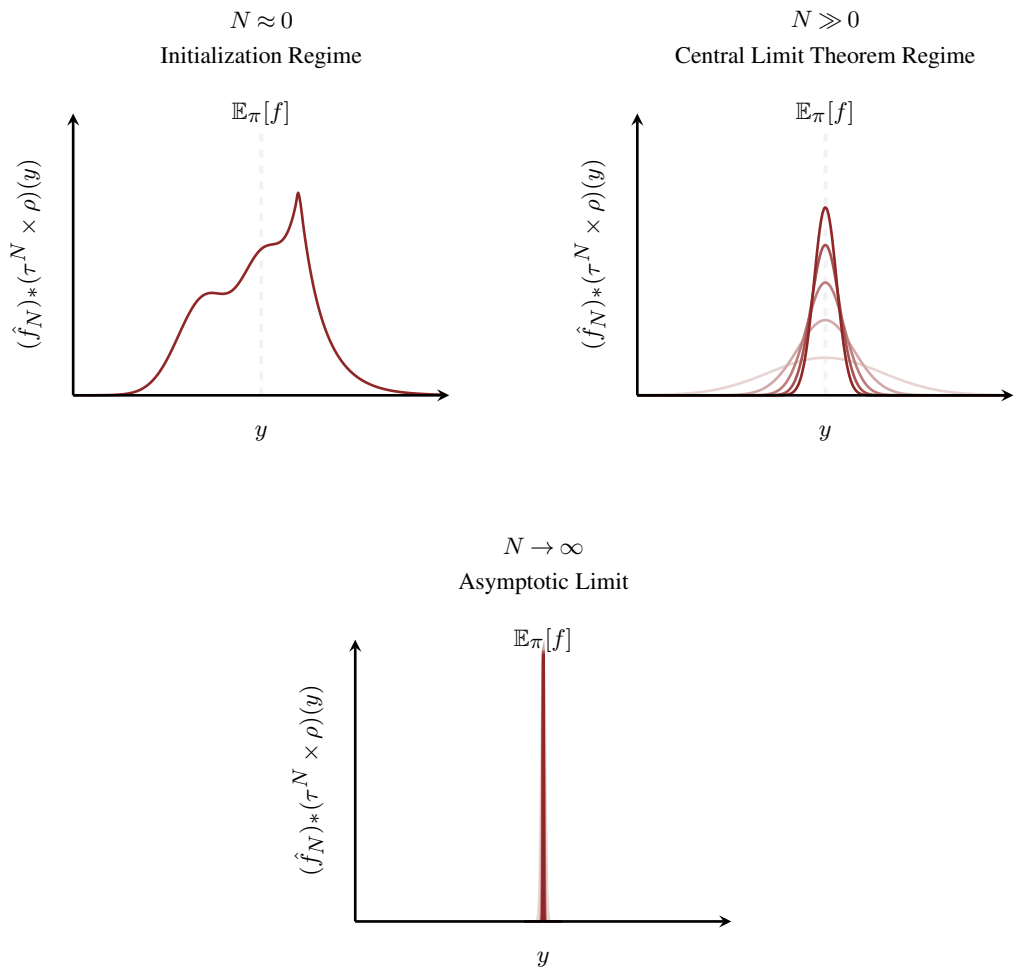

For small $N$ all of these estimates suffer from a bias that monotonically decays with $N^{-1}$.
As $N$ increases the bias becomes negligible and the Markov chain Monte Carlo central limit 
theorem kicks in, allowing us to empirically quantify the error of the estimator using only 
the history of the realized Markov chain.  Letting $N$ grow to infinity the normal approximation 
given by the central limit theorem continues to narrow until it finally converges to a Dirac 
distribution in the asymptotic limit (Figure \ref{fig:phases_of_convergence}).

\section{Acknowledgements}

I thank Sam Livingstone, Simon Byrne, Sam Power, and Finn Lindgren for helpful comments 
while taking all responsibility for any errors and unconventional mathematical notations. 

\bibliography{mcmc_convergence}

\begin{thebibliography}{13}

\bibitem[\protect\citeauthoryear{Chan and Geyer}{1994}]{ChanEtAl:1994}
\begin{barticle}[author]
\bauthor{\bsnm{Chan},~\bfnm{Kung~Sik}\binits{K.~S.}} \AND
  \bauthor{\bsnm{Geyer},~\bfnm{Charles~J.}\binits{C.~J.}}
(\byear{1994}).
\btitle{Discussion: Markov Chains for Exploring Posterior Distributions}.
\bjournal{The Annals of Statistics}
\bvolume{22}
\bpages{1747--1758}.
\end{barticle}
\endbibitem

\bibitem[\protect\citeauthoryear{Cover and Thomas}{2006}]{CoverEtAl:2006}
\begin{bbook}[author]
\bauthor{\bsnm{Cover},~\bfnm{Thomas~M.}\binits{T.~M.}} \AND
  \bauthor{\bsnm{Thomas},~\bfnm{Joy~A.}\binits{J.~A.}}
(\byear{2006}).
\btitle{Elements of information theory},
\bedition{Second} ed.
\bpublisher{Wiley-Interscience [John Wiley \& Sons], Hoboken, NJ}.
\end{bbook}
\endbibitem

\bibitem[\protect\citeauthoryear{Diaconis and
  Freedman}{1999}]{DiaconisEtAl:1999}
\begin{barticle}[author]
\bauthor{\bsnm{Diaconis},~\bfnm{Persi}\binits{P.}} \AND
  \bauthor{\bsnm{Freedman},~\bfnm{David}\binits{D.}}
(\byear{1999}).
\btitle{Iterated Random Functions}.
\bjournal{SIAM review}
\bvolume{41}
\bpages{45--76}.
\end{barticle}
\endbibitem

\bibitem[\protect\citeauthoryear{Dudley}{2002}]{Dudley:2002}
\begin{bbook}[author]
\bauthor{\bsnm{Dudley},~\bfnm{R.~M.}\binits{R.~M.}}
(\byear{2002}).
\btitle{Real analysis and probability}.
\bseries{Cambridge Studies in Advanced Mathematics}
\bvolume{74}.
\bpublisher{Cambridge University Press, Cambridge}.
\end{bbook}
\endbibitem

\bibitem[\protect\citeauthoryear{Harris}{1956}]{Harris:1956}
\begin{binproceedings}[author]
\bauthor{\bsnm{Harris},~\bfnm{T.~E.}\binits{T.~E.}}
(\byear{1956}).
\btitle{The existence of stationary measures for certain {M}arkov processes}.
In \bbooktitle{Proceedings of the {T}hird {B}erkeley {S}ymposium on
  {M}athematical {S}tatistics and {P}robability, 1954--1955, vol. {II}}
\bpages{113--124}.
\bpublisher{University of California Press, Berkeley and Los Angeles}.
\end{binproceedings}
\endbibitem

\bibitem[\protect\citeauthoryear{Joulin and Ollivier}{2010}]{JoulinEtAl:2010}
\begin{barticle}[author]
\bauthor{\bsnm{Joulin},~\bfnm{Ald{\'e}ric}\binits{A.}} \AND
  \bauthor{\bsnm{Ollivier},~\bfnm{Yann}\binits{Y.}}
(\byear{2010}).
\btitle{Curvature, Concentration and Error Estimates for {M}arkov {C}hain
  {M}onte {C}arlo}.
\bjournal{The Annals of Probability}
\bvolume{38}
\bpages{2418--2442}.
\end{barticle}
\endbibitem

\bibitem[\protect\citeauthoryear{Kantorovi\v{c} and
  Rubin\v{s}te\u{\i}n}{1958}]{KantorovicEtAl:1958}
\begin{barticle}[author]
\bauthor{\bsnm{Kantorovi\v{c}},~\bfnm{L.~V.}\binits{L.~V.}} \AND
  \bauthor{\bsnm{Rubin\v{s}te\u{\i}n},~\bfnm{G.~\v{S}.}\binits{G.~v.}}
(\byear{1958}).
\btitle{On a space of completely additive functions}.
\bjournal{Vestnik Leningrad. Univ.}
\bvolume{13}
\bpages{52--59}.
\end{barticle}
\endbibitem

\bibitem[\protect\citeauthoryear{Lee}{2018}]{Lee:2018}
\begin{bbook}[author]
\bauthor{\bsnm{Lee},~\bfnm{John~M.}\binits{J.~M.}}
(\byear{2018}).
\btitle{Introduction to {R}iemannian manifolds}.
\bseries{Graduate Texts in Mathematics}
\bvolume{176}.
\bpublisher{Springer, Cham}.
\end{bbook}
\endbibitem

\bibitem[\protect\citeauthoryear{M\"{u}ller}{1997}]{Muller:1997}
\begin{barticle}[author]
\bauthor{\bsnm{M\"{u}ller},~\bfnm{Alfred}\binits{A.}}
(\byear{1997}).
\btitle{Integral probability metrics and their generating classes of
  functions}.
\bjournal{Adv. in Appl. Probab.}
\bvolume{29}
\bpages{429--443}.
\end{barticle}
\endbibitem

\bibitem[\protect\citeauthoryear{Ollivier}{2009}]{Ollivier:2009}
\begin{barticle}[author]
\bauthor{\bsnm{Ollivier},~\bfnm{Yann}\binits{Y.}}
(\byear{2009}).
\btitle{{R}icci Curvature of {M}arkov Chains on Metric Spaces}.
\bjournal{Journal of Functional Analysis}
\bvolume{256}
\bpages{810--864}.
\end{barticle}
\endbibitem

\bibitem[\protect\citeauthoryear{Roberts and
  Rosenthal}{2004}]{RobertsEtAl:2004}
\begin{barticle}[author]
\bauthor{\bsnm{Roberts},~\bfnm{Gareth~O}\binits{G.~O.}} \AND
  \bauthor{\bsnm{Rosenthal},~\bfnm{Jeffrey~S}\binits{J.~S.}}
(\byear{2004}).
\btitle{General State Space {M}arkov Chains and {MCMC} Algorithms}.
\bjournal{Probability Surveys}
\bvolume{1}
\bpages{20--71}.
\end{barticle}
\endbibitem

\bibitem[\protect\citeauthoryear{Sriperumbudur
  et~al.}{2009}]{SriperumbudurEtAl:2009}
\begin{barticle}[author]
\bauthor{\bsnm{Sriperumbudur},~\bfnm{Bharath~K.}\binits{B.~K.}},
  \bauthor{\bsnm{Fukumizu},~\bfnm{Kenji}\binits{K.}},
  \bauthor{\bsnm{Gretton},~\bfnm{Arthur}\binits{A.}},
  \bauthor{\bsnm{Sch{\"o}lkopf},~\bfnm{Bernhard}\binits{B.}} \AND
  \bauthor{\bsnm{Lanckriet},~\bfnm{Gert R.~G.}\binits{G.~R.~G.}}
(\byear{2009}).
\btitle{On integral probability metrics, $\phi$-divergences and binary
  classification}.
\bjournal{arXiv e-prints}
\bvolume{0901.2698}.
\end{barticle}
\endbibitem

\bibitem[\protect\citeauthoryear{Tierney}{1994}]{Tierney:1994}
\begin{barticle}[author]
\bauthor{\bsnm{Tierney},~\bfnm{Luke}\binits{L.}}
(\byear{1994}).
\btitle{Markov chains for exploring posterior distributions}.
\bjournal{Ann. Statist.}
\bvolume{22}
\bpages{1701--1762}.
\end{barticle}
\endbibitem

\end{thebibliography}
\bibliographystyle{imsart-nameyear}

\end{document}